\documentclass[12pt]{amsart}
\usepackage{amsmath,amscd,amssymb,amsfonts}
\newcommand{\h}{\hbox}
\newcommand{\q}{\quad}
\newcommand{\nin}{\par\noindent}
\newcommand{\bs}{\par\bigskip}
\newcommand{\ms}{\par\medskip}
\newcommand{\sk}{\par\smallskip} 
\newcommand{\ssb}{\raise.2ex\h{${\scriptscriptstyle\bullet}$}}
\newcommand{\ssc}{\,\raise.2ex\hbox{${\scriptstyle\circ}$}\,}
\newcommand{\msum}{\h{$\sum$}}
\newcommand{\mprod}{\h{$\prod$}}
\newcommand{\mopl}{\h{$\bigoplus$}}
\newcommand{\mcup}{\h{$\bigcup$}}

\newcommand{\mwdg}{\h{$\bigwedge$}}
\newcommand{\be}{{\mathbf e}}
\newcommand{\C}{{\mathbb C}}
\newcommand{\N}{{\mathbb N}}
\newcommand{\bP}{{\mathbb P}}
\newcommand{\Q}{{\mathbb Q}}
\newcommand{\R}{{\mathbb R}}

\newcommand{\Z}{{\mathbb Z}}

\newcommand{\Ac}{{\mathcal A}}
\newcommand{\B}{{\mathcal B}}
\newcommand{\D}{{\mathcal D}}
\newcommand{\Ec}{{\mathcal E}}
\newcommand{\Gc}{{\mathcal G}}
\newcommand{\Ic}{{\mathcal I}}
\newcommand{\Jc}{{\mathcal J}}
\newcommand{\Lc}{{\mathcal L}}
\newcommand{\M}{{\mathcal M}}
\newcommand{\Oc}{{\mathcal O}}

\newcommand{\alt}{\widetilde{\alpha}}
\newcommand{\tb}{\widetilde{b}}
\newcommand{\tD}{\widetilde{D}}
\newcommand{\tE}{\widetilde{E}}
\newcommand{\tf}{\widetilde{f}}
\newcommand{\tH}{\widetilde{H}}

\newcommand{\tm}{\widetilde{m}}
\newcommand{\tP}{\widetilde{P}}
\newcommand{\tR}{\widetilde{R}}
\newcommand{\tX}{\widetilde{X}}
\newcommand{\tY}{\widetilde{Y}}
\newcommand{\tZ}{\widetilde{Z}}
\newcommand{\tcG}{\widetilde{\mathcal G}}
\newcommand{\tal}{\widetilde{\alpha}}

\newcommand{\red}{\text{{\rm red}}}

\newcommand{\Gr}{\text{{\rm Gr}}}
\newcommand{\Ker}{\text{{\rm Ker}}}
\newcommand{\DR}{\text{{\rm DR}}}

\newcommand{\GCD}{\h{\rm GCD}}
\newcommand{\JC}{\h{\rm JC}}
\newcommand{\mult}{\h{\rm mult}}
\newcommand{\Sing}{\text{{\rm Sing}}\,}
\newcommand{\supp}{\text{{\rm supp}}\,}
\newcommand{\simto}{\buildrel\sim\over\longrightarrow}
\newcommand{\al}{\alpha}
\newcommand{\la}{\lambda}
\newcommand{\lab}{\bar{\lambda}}
\newcommand{\om}{\omega}
\newcommand{\Si}{\Sigma}
\newcommand{\dd}{\partial}
\newcommand{\ddd}{{\rm d}}

\newcommand{\kod}{\h{$\frac{k}{d}$}}
\newcommand{\wdg}{\,{\wedge}\,}
\newcommand{\bl}{\bigl}
\newcommand{\br}{\bigr}
\newcommand{\into}{\hookrightarrow}
\newcommand{\les}{\leqslant}
\newcommand{\ges}{\geqslant}
\begin{document}
\title[Bernstein-Sato Polynomials]
{Bernstein-Sato Polynomials of Hyperplane Arrangements}
\author[M. Saito]{Morihiko Saito}
\address{RIMS Kyoto University, Kyoto 606-8502 Japan}
\dedicatory{Dedicated to Joseph Bernstein}
\begin{abstract} We show that the Bernstein-Sato polynomial (that is, the $b$-function) of a hyperplane arrangement with a reduced equation is calculable by combining a generalization of Malgrange's formula with the theory of Aomoto complexes due to Esnault, Schechtman, Terao, Varchenko, and Viehweg in certain cases. We prove in general that the roots are greater than $-2$ and the multiplicity of the root $-1$ is equal to the (effective) dimension of the ambient space. We also give an estimate of the multiplicities of the roots in terms of the multiplicities of the arrangement at the dense edges, and provide a method to calculate the Bernstein-Sato polynomial at least in the case of 3 variables with degree at most 7 and generic multiplicities at most 3. Using our argument, we can terminate the proof of a conjecture of Denef and Loeser on the relation between the topological zeta function and the Bernstein-Sato polynomial of a reduced hyperplane arrangement in the 3 variable case.
\end{abstract}
\maketitle
\centerline{\bf Introduction}
\bs\nin
Let $D$ be a hypersurface of a smooth affine algebraic variety $X$ with a defining equation $f$. The Bernstein-Sato polynomial (that is, the $b$-function) $b_f(s)$ is the monic polynomial of the lowest degree satisfying
$$b_f(s)f^{s}=Pf^{s+1}\q\h{for}\,\,\,P\in\D_{X}[s],$$
where $\D_{X}$ is the ring of differential operators, see \cite{Be}, \cite{SatSh}, etc.
Let $R_f$ be the set of roots of $b_f(-s)$, and $m_{\al}$ be the multiplicity of $\al\in R_f$. Set $\al_f:=\min R_f$. This coincides with the minimal jumping coefficient, see \cite{ELSV}, \cite{Ko} (and also \cite{Sa4} for the real case, where the analytic continuation as in \cite{Be} is used in an essential way). Let $n=\dim X$. Then $R_f\subset\Q_{>0}$ (see \cite{Ka1}), and $m_{\al}\les n$ since $b_f(s)$ is closely related to the monodromy on the nearby cycle sheaf $\psi_f\C_{X}$, see \cite{Ka2}, \cite{Ma2}. Define $\tR_f,\tm_{\al},\tal_f$ by replacing $b_f(s)$ with the reduced (or microlocal) Bernstein-Sato polynomial $\tb_f(s):=b_f(s)/(s+1)$. Then we have more precisely (see \cite{Sa2}):
$$\tR_f\subset [\tal_f,n-\tal_f],\q\tm_{\al}\les n-\tal_f-\al+1.
\leqno(0.1)$$
It is well-known that the first inclusion is optimal in the weighted homogeneous isolated singularity case where $\max\tR_f=n-\tal_f$ and $\tm_{\al}\les 1$ (see, for instance, \cite[Chapters~1 and 4]{Sat} (in Japanese) or \cite{Ma1}). In fact, in the case $D$ has an isolated singularity more generally, it was shown by Malgrange (\h{\it loc.~cit.}) that $\tb_f(s)$ coincides with the minimal polynomial of $-\dd_tt$ on $\tH''_{\!f,0}/t\tH''_{\!f,0}$ where $\tH''_{\!f,0}$ is the saturation of the Brieskorn lattice $H''_{\!f,0}$ of $f$ at $0$ (see \cite{Br1}). In the non-isolated singularity case, however, it is quite difficult to calculate explicitly the Bernstein-Sato polynomial without using computer programs.
\sk
Assume now $D$ is a central hyperplane arrangement in $X=\C^n$, that is, the affine cone of a projective hyperplane arrangement $Z$ in $Y=\bP^{n-1}$. We may assume that $D$ is not the pull-back of an arrangement in a vector space of strictly lower dimension. (In this case, $n$ is called the {\it effective} dimension.) Set $d:=\deg f$. Then
\ms\nin
{\bf Theorem~1.}\,\, $\max R_f<2-\frac{1}{d},\,\,\, m_1=n$.
\ms
This is quite different from (0.1) in the case of weighted homogeneous isolated singularities explained above. Theorem~1 implies that $1$ is the only integral root of $b_f(-s)$ for hyperplane arrangements as is shown in \cite{Wa1}. For the proof of Theorem~1, we use a generalization of the above formula of Malgrange~\cite{Ma1} (see (1.4) below) together with the theory of Aomoto complexes due to Esnault, Schechtman, Terao, Varchenko, Viehweg (\cite{ESV}, \cite{STV}) for the calculation of certain twisted de Rham cohomology groups, see (2.2) below. Note that the latter enables us to reduce a problem of linear algebra for infinite dimensional vector spaces to the one for finite dimensional vector spaces, see Remark~(1.3) below.
\sk
As for $\al_f=\min R_f$, we have the following by a theorem of Musta\c{t}\v{a} on the jumping coefficients \cite{Mu} combined with \cite{Ko}, \cite{ELSV}:
$$\al_f=\min\bl(\h{$\frac{n}{d}$},\al'_f\br),
\leqno(0.2)$$
where $\al'_f=\min\bigcup_{x\ne 0}R_{f,x}$ with $R_{f,x}$ the set of roots of the local Bernstein-Sato polynomial $b_{f,x}(s)$ of $f$ at $x$ up to a sign, see (1.8) below.
\sk
In order to describe the non-integral roots of the Bernstein-Sato polynomial, we need some terminology from \cite{STV}. We say that an edge $L$ of $ D$ (that is, an intersection of irreducible components of $D$) is {\it dense} if the irreducible components $D_i$ containing $L$ are identified with an indecomposable arrangement. (See (3.4) below for decomposable arrangements.) Let $\D\Ec(D)$ be the set of dense edges of $D$. Let $m_{L}$ be the number of $D_i$ containing $L$. For $\la\in\C$, let $\D\Ec(D,\la)$ be the subset of $\D\Ec(D)$ consisting of $L$ such that $\la^{m_{L}}=1$. We say that $L, L'\in\D\Ec(D,\la)$ are {\it strongly adjacent} if $L\subset L'$ or $L\supset L'$ or $L\cap L'$ is non-dense. Let $m(\la)$ be the maximal number of the elements of subsets $S$ of $\D\Ec(D,\la)$ such that any two elements of $S$ are strongly adjacent. Using an embedded resolution of singularities in \cite{STV} together with \cite{Ka2}, \cite{Ma2}, we get
\ms\vbox{\nin
{\bf Theorem~2.} {\it We have $m_{\al}\les m(\la)$ with $\la=\exp(-2\pi i\al)$. In particular, we have the inclusion
$$R_f\subset\bigcup_{L\in\D\Ec(D)}\frac{\Z}{m_L},$$
and $m_{\al}=1$ for any $\al\in R_f\setminus\Z$ if $\GCD(m_{L},m_{L'})=1$ for any dense edges $L, L'$ of $D$ such that $L\subset L'$ or $L\cap L'$ is non-dense.}}
\ms
If $f$ is generic (that is, if $Z$ is a divisor with normal crossings) and $d>n\ges 2$, then U.~Walther \cite{Wa1} proved (except for the multiplicity of $-1$):
$$b_f(s)=(s+1)^n\,\prod_{k=n}^{d-1}\bl(s+\kod\br)\prod_{k=d+1}^{2d-2}\bl(s+\kod\br).
\leqno(0.3)$$
Note that Theorems~1--2 and (0.2) imply that $b_f(s)$ divides the right-hand side since the dense edges in this case consist of the irreducible components $D_i$ and the origin $0$ where $m_{D_i}=1$ and $m_{0}=d$. We can prove the equality using a calculation of the spectrum as in \cite{Sa3}. We also calculate some examples in the non-generic case, see Theorems~(4.9) and (4.11--12), and Examples~(4.13) below. Note that there is no combinatorial formula for the Bernstein-Sato polynomials of non-generic hyperplane arrangements (see \cite{Wa2}).
\sk
Finally our arguments can be used to study the following conjecture of J.~Denef and F.~Loeser on the relation between the Bernstein-Sato polynomial $b_f(s)$ and the topological zeta function $Z^{\rm top}_{f,0}(s)$ (see \cite{DeLo}):
\ms\nin
{\bf Conjecture~1.} The product of the Bernstein-Sato polynomial and the topological zeta function $b_f(s)\,Z^{\rm top}_{f,0}(s)$ has no poles.
\ms
In fact, we have the following.
\ms\nin
{\bf Theorem~3.} {\it Conjecture~$1$ is true for reduced hyperplane arrangements with $n=3$.}
\ms
Note that Theorem~3 was proved in \cite{BSY} with $b_f(s)\,Z^{\rm top}_{f,0}(s)$ replaced by $b_f(s)^k\,Z^{\rm top}_{f,0}(s)$ for $k\gg 0$. (This means that any pole of $Z^{\rm top}_{f,0}(s)$ is a root of $b_f(s)$ with multiplicity forgotten.) Moreover the proof of Theorem~3 was reduced by \cite[Proposition~2.4]{BSY} to the following (see Remark~(5.4)(i) below):
\ms\nin
{\bf Proposition~1.} {\it Assume that $d=3m$ with $m\in\N$, and there is a point of $Z\subset\bP^2$ with multiplicity $2m$. Then $-\frac{1}{m}$ is a root of $b_f(s)$ with multiplicity $2$.}
\ms
We can prove this by using Theorem~(3.8) below together with the theory of Aomoto complexes (see (5.1) below).
\sk
We thank N.~Budur, A.~Dimca, and M.~Musta\c{t}\v{a} for useful discussions and comments about the subject of this paper. We are grateful to Professor M.~Noro who has calculated some examples as in Examples~(4.13) using his computer program Risa/Asir, and verified the coincidence with the calculation in this paper. We also thank A.~Dimca and G.~Sticlaru for their calculation of Hilbert series in many examples by using a Singular computer program, which also verifies the correctness of our arguments.
This work was supported by JSPS Kakenhi 17540023 and 24540039.
\sk
In Section~1 we review some basic facts from the theory of multiplier ideals, Bernstein-Sato polynomial and spectrum, see \cite{Sa3}. In Section~2 we review and partially generalize the theory of Aomoto complexes due to Esnault, Schechtman, Terao, Varchenko, Viehweg. This can be combined effectively with the construction in (1.5). In Section~3 we prove the main theorems, and calculate the nearby and vanishing cycle sheaves in certain cases. In Section~4 we show how to calculate the Bernstein-Sato polynomial in certain cases including the generic one. In Section~5 we terminate the proof of Theorem~3 by showing Proposition~1.
\bs\bs
\vbox{\centerline{\bf 1. Preliminaries}
\bs\nin
In this section we review some basic facts from the theory of multiplier ideals, Bernstein-Sato polynomial and spectrum, see \cite{Sa3}.}
\ms\nin
{\bf 1.1.~Bernstein-Sato polynomials.} Let $X$ be a complex manifold, and $D$ be a hypersurface defined by a holomorphic function $f$. We have a canonical injection of $\D_{X}[s]$-modules
$$M:=\D_{X}[s]f^{s}\subset\B_f:=\Oc_{X}\otimes_{\C}\C[\dd_t],$$
such that $f^{s}$ is sent to $1\otimes 1$, where $s=-\dd_tt$, see \cite{Ka1}, \cite{Ma1}. Note that $\B_f$ is the direct image of $\Oc_{X}$ by the graph embedding $i_f:X\to X\times\C$ as a $\D$-module, and the action of $\D_{X\times\C}$ on $\B_f$ is defined by identifying $1\otimes 1$ with the delta function $\delta(t-f)$.
\sk
The Bernstein-Sato polynomial (that is, the $b$-function) $b_f(s)$ is the minimal polynomial of the action of $s$ on $M/tM$. Since $M/tM$ is holonomic, the Bernstein-Sato polynomial exists if $X$ is (relatively) compact or if $X, f$ are algebraic. This implies also that the Bernstein-Sato polynomial of an algebraic function coincides with the Bernstein-Sato polynomial of the associated analytic function. Restricting to the stalk at $x\in D$, we can also define the local Bernstein-Sato polynomial $b_{f,x}(s)$. Note that the global Bernstein-Sato polynomial is the least common multiple of the local Bernstein-Sato polynomials if $X$ is affine or Stein. In the case $D$ is an affine cone, the global Bernstein-Sato polynomial coincides with the local Bernstein-Sato polynomial at $0$.
\sk
By Kashiwara \cite{Ka2} and Malgrange \cite{Ma2}, $\B_f$ has the filtration $V$ (indexed by $\Q$) together with a canonical isomorphism of perverse sheaves
$$\DR_{X}(\Gr_V^{\al}\B_f)=\psi_{f,\la}\C_{X}[n{-}1]\q\h{for}\,\,\,\al>0,\,\la=\exp(-2\pi i\al),
\leqno(1.1.1)$$
such that $\exp(-2\pi i\dd_tt)$ on the left-hand side corresponds to the monodromy $T$ on the right-hand side. Here $\DR_{X}$ is the de Rham functor inducing an equivalence of categories between regular holonomic $\D$-modules and perverse sheaves, and $\psi_{f,\la}\C_{X}=\Ker(T_{s}-\la)\subset\psi_f\C_{X}$ in the abelian category of shifted perverse sheaves \cite{BBD}, where $\psi_f\C_{X}$ is the nearby cycle sheaf \cite{De3}, and $T_{s}$ is the semisimple part of the monodromy.
\sk
Let $G$ be the increasing filtration on $\Gr_V^{\al}\B_f\, (\al\in (0,1])$ defined by
$$t^i(G_i\Gr_V^{\al}\B_f)=\Gr_V^{\al+i}M\subset\Gr_V^{\al+i}\B_f\q\bl(i\in\N\br),
\leqno(1.1.2)$$
where $G_i\Gr_V^{\al}\B_f=0$ for $i<0$ (since $M\subset V^{>0}\B_f$ by \cite{Ka1}). Then
$$t^i:\Gr_i^{G}\Gr_V^{\al}\B_f\simto\Gr_V^{\al+i}(M/tM)\q\bl(i\in\N\br),$$
and we get
$$m_{\al+i}=\min\{k\in\N\,|\, N^{k}\Gr^{G}_i\Gr_V^{\al}\B_f=0\}\q\bl(\al\in(0,1], i\in\N\br),
\leqno(1.1.3)$$
where $N=-(\dd_tt-\al)$. In particular, $m_{\al}\les n$ for any $\al$.
\ms\nin
{\bf 1.2.~Hodge and pole order filtrations.} Let $X, D, f$ be as in (1.1), and $F_{\!f,x}$ be the Milnor fiber of $f$ at $x\in D$, that is,
$$F_{\!f,x}=\{z\in X\,\big|\,|z-x|<\varepsilon,f(z)=\delta\}\q\h{for}\,\,\, 0<\delta\ll\varepsilon\ll 1,$$
where $|z-x|$ is defined by choosing a local coordinate system of $X$. Then it is well-known that the Milnor cohomology $H^j(F_{\!f,x},\C)$ has a canonical mixed Hodge structure. This can be defined, for instance, by using \cite{Sa1} since we have a canonical isomorphism
$$R^ji_{x}^*\psi_f\C_{X}=H^j(F_{\!f,x},\C),$$
where $i_{x}:\{x\}\to X$ is the natural inclusion. Let $T_{s}$ be the semisimple part of the monodromy $T$. We have the decompositions
$$\psi_f\C_{X}=\mopl_{\la}\psi_{f,\la}\C_{X},\q H^j(F_{\!f,x},\C)=\mopl_{\la}H^j(F_{\!f,x},\C)_{\la},$$
such that the action of $T_{s}$ on $\psi_{f,\la}\C_{X}, H^j(F_{\!f,x},\C)_{\la}$ is the multiplication by $\la\in\C^*$.
\sk
Let $\Gc_f$ be the Gauss-Manin system of the highest degree associated with $f$ at $x\in D$. Let $\Gc_f^{(0)}\subset\Gc_f$ be the Brieskorn lattice \cite{Br1}, that is,
$$\Gc_f^{(0)}=H''_{f,x}/(t\h{-torsion})\q\h{with}\q H''_{\!f,x}:=\Omega_{X,x}^n/\ddd\Ac^{n-1}_{X,x},$$
where $\Ac^j_{X,x}:={\rm Ker}\bl(\ddd f\wdg:\Omega_{X,x}^j\to\Omega_{X,x}^{j+1}\br)$.
Note that $\Gc_f$ is the localization of $\Gc^{(0)}$ by the action of $\dd_t^{-1}$. Here $\dd_t^{-1}[\om]$ is defined by $[\ddd f\wdg\eta]$ with $\ddd\eta=\om$, and the action of $t$ is defined by the multiplication by $f$. It is known that the $t$-torsion of $H''_f$ coincides with the $\dd_t^{-1}$-torsion, see, for instance, \cite{BaSa}. Let $\tcG_f^{(0)}$ be the saturation of $\Gc_f^{(0)}$, that is,
$$\tcG_f^{(0)}=\msum_{i\ges 0}\,(\dd_tt)^i\Gc_f^{(0)}=\msum_{i\ges 0}\,(t\dd_t)^i\Gc_f^{(0)}.$$
Set
$$\tcG_f^{(-i)}=\dd_t^i\tcG_f^{(0)},\q\Gc_f^{(-i)}=\dd_t^i\Gc_f^{(0)}.$$
\sk
Let $V$ be the filtration of Kashiwara and Malgrange. We have canonical isomorphisms
$$H^{n-1}(F_{\!f,x},\C)_{\la}=\Gr_V^{\al}\Gc_f\q\h{for}\,\,\la=\exp(-2\pi i\al)\,\,\h{with}\,\,\al\in (0,1],
\leqno(1.2.1)$$
such that the monodromy $T$ corresponds to $\exp(-2\pi i\dd_tt)$. Using the (1.2.1), we define decreasing filtrations $P$ and $\tP$ on the Milnor cohomology so that
$$\aligned P^{n-1-i} H^{n-1}(F_{\!f,x},\C)_{\la} &=\Gr_V^{\al}\Gc_f^{(-i)}\,(\subset\Gr_V^{\al}\Gc_f),\\
\tP^{n-1-i} H^{n-1}(F_{\!f,x},\C)_{\la} &=\Gr_V^{\al}\tcG_f^{(-i)}\,(\subset\Gr_V^{\al}\Gc_f).\endaligned
\leqno(1.2.2)$$
Note that, if there is a vector field $\xi$ such that $\xi f=f$, then
$$\tcG_f^{(-i)}=\Gc_f^{(-i)},\q\tP=P.
\leqno(1.2.3)$$
\ms\nin
{\bf Remark~1.3.} In general it is quite difficult to calculate the filtration $P$ even if $D$ is an affine cone. In this case we have the decomposition of the Brieskorn module at the origin $H''_{\!f,0}=\widehat{\bigoplus}_{i\in\N}(H''_{\!f,0})_i$ by using the degree such that $\deg x_i=\deg dx_i=1$. Then the action of $\dd_tt$ on $(H''_{\!f,0})_i$ is the multiplication by $i/d$ where $d=\deg f$, and each $(H''_{\!f,0})_i$ can be calculated by using only finite dimensional vector spaces. However, we have to consider the inductive system $\{(H''_{\!f,0})_{i+kd}\}_{k\in\N}$ defined by the multiplication by $f$ to kill the $t$-torsion, and we get a problem concerning infinite dimensional vector spaces.
\ms
The following is a generalization of \cite{Ma1}, and is proved in \cite[Theorem~2]{Sa3}.
\ms\vbox{\nin
{\bf Theorem~1.4.} {\it The filtration $\tP$ on $H^{n-1}(F_{\!f,x},\C)_{\la}$ contains the Hodge filtration $F$, and for any $\al\in\Q$ such that $\la=\exp(-2\pi i\al)$, we have the following\,{\rm :}
\sk\nin
$(i)$ If $\Gr_{\tP}^{p}H^{n-1}(F_{\!f,x},\C)_{\la}\ne 0$ with $p=[n-\al]$, then $\al$ is a root of $b_{f,x}(-s)$.
\sk\nin
$(ii)$ If $\al+i$ is not a root of $b_{f,y}(-s)$ for any $y\ne x$ and any $i\in\N$, then the converse of the assertion $(i)$ holds.
\sk\nin
$(iii)$ If $\la$ is not an eigenvalue of the Milnor monodromy at $y\ne x$, then the multiplicity of the root $\al$ coincides with the degree of the minimal polynomial of the action of the monodromy on $\Gr_{\tP}^{p}H^{n-1}(F_{\!f,x},\C)_{\la}$.}}
\ms\nin
{\bf 1.5.~Affine cone case.} Assume $X=\C^n$ with $n\ges 3$, and $D$ is the affine cone of a divisor $Z$ of degree $d$ on $Y:=\bP^{n-1}$. By \cite[Ch.~6, Thm.~2.9]{Di} (or \cite[Section 1.8]{DiSa2}), the pole order filtration $P$ in (1.2) can be calculated in this case as below. Here we have $\tP=P$, since $f$ is a homogeneous polynomial, see (1.2.3).
\sk
We have a natural cyclic covering $\pi:\tY\to Y$ of degree $d$ which is ramified along $Z$, and such that its restriction over $U:=Y\setminus Z$ is isomorphic to the restriction of the natural projection $\C^n\setminus\{0\}\to\bP^{n-1}$ to the Milnor fiber $F_{\!f,0}:=f^{-1}(1)$ of a polynomial $f$ defining the affine cone $D$ of $Z$. Here the geometric Milnor monodromy corresponds to a generator of the covering transformation group of $\pi:F_{\!f,0}\to U$ (see also \cite[Section 1.3]{BuSa2}).
\sk
For $k=1,\dots, d$, let $L^{(k/d)}$ be the direct factor of $\pi_*\C_{F_{\!f,0}}$ on which the action of the Milnor monodromy is the multiplication by $\exp(-2\pi ik/d)$. Then $L^{(k/d)}$ is a local system of rank 1 on $U$, and
$$H^j(U,L^{(k/d)})=H^j(F_{\!f,0},\C)_{\la},
\leqno(1.5.1)$$
where $\la=\exp(-2\pi ik/d)$. Let $\Lc^{(k/d)}$ be the meromorphic extension of $L^{(k/d)}\otimes_{\C}\Oc_{U}$. This is a regular holonomic $\D_{Y}$-module on which the action of a function $h$ defining $Z$ is bijective. We see that $\Lc^{(k/d)}$ is locally isomorphic to a free $\Oc_{Y}(*Z)$-module generated by a multivalued function $h_j^{-k/d}$ where $h_j=x_j^{-d}f$ on $\{x_j\ne 0\}\subset\bP^{n-1}$. Note that the $\Oc_{Y}$-submodule generated locally by $h_j^{-k/d}$ is isomorphic to $\Oc_{Y}(k)$.
\sk
The pole order filtration $P_i\Lc^{(k/d)}$ is defined to be the locally free $\Oc_{Y}$-submodule of $\Lc^{(k/d)}$ generated by $h_j^{-i-(k/d)}$ on $\{x_j\ne 0\}$ for $i\in\N$, and $P_i\Lc^{(k/d)}=0$ for $i<0$. We have
$$P_i\Lc^{(k/d)}\cong\Oc_{Y}(id+k)\q\h{for}\,\,\,\,i\in\N,\,\,k\in\{1,\dots,d\}.
\leqno(1.5.2)$$
On the other hand, there is the Hodge filtration $F$ on $\Lc^{(k/d)}$ such that $F_i\Lc^{(k/d)}=P_i\Lc^{(k/d)}$ outside $\Sing Z_{\red}$ for any $i$ by the theory of mixed Hodge modules. Then we have $F_i\Lc^{(k/d)}\subset P_i\Lc^{(k/d)}$ on $Y$ since $P_i\Lc^{(k/d)}$ is locally free and $\Sing Z_{\red}$ has codimension $\ges 2$ in $Y$, see also \cite[Section 1.8]{DiSa2}.
\sk
The Hodge and pole order filtrations are closely related respectively to the spectrum \cite{St3} and the Bernstein-Sato polynomial of $f$. Indeed, the Hodge filtration $F$ on $\Lc^{(k/d)}$ induces the Hodge filtration on the Milnor cohomology by taking the de Rham cohomology. Similarly the pole order filtration $P$ on the Milnor cohomology is defined by using the de Rham cohomology. Here the filtration is shifted by the degree of the differential forms $n-1$, and the associated decreasing filtration is used.
\ms
The following two propositions are proved in \cite{Sa3}.
\ms\vbox{\nin
{\bf Proposition~1.6.} {\it Under the assumption of $(1.5)$ the pole order filtration $P$ in {\rm (1.5)} coincides with the filtration $\tP=P$, see {\rm (1.2.3)}. Moreover, for $\al=k/d\in (0,1)$ and $\la=\exp(-2\pi i\al)$, we can identify $P^{n-1-j}$ for $j\in\N$ in {\rm (1.5)} with the image of
$$\Gr_V^{\al+j}\Gc_f^{(0)}\subset\Gr_V^{\al+j}\Gc_f\simeq\Gr_V^{\al}\Gc_f\simeq H^{n-1}(F_{\!f,0},\C)_{\la},$$
where the middle isomorphism can be induced by both $\dd_t^j$ and $t^{-j}$, and the last morphism is induced by {\rm (1.2.1)}.}}
\ms\vbox{\nin
{\bf Proposition~1.7.} {\it With the notation of {\rm (1.2)}, assume $\Gr_{n-1+k}^WH^{n-1}(F_{\!f,x},\C)_{\la}\ne 0$ for a positive integer $k$, where $W$ is the weight filtration. Then $N^{k}\ne 0$ on $\psi_{f,\la}\C_{X}$ in the category of shifted perverse sheaves, where $N$ is the logarithm of the unipotent part of the monodromy $T$.}}
\ms\nin
{\bf 1.8.~Multiplier ideals and spectrum.} Let $D\subset X$ be the affine cone of a divisor $Z\subset\bP^{n-1}$ with a defining equation $f$, where $n\ges 3$. Let $R_{f,x}$ be the set of the roots of $b_{f,x}(-s)$, and set
$$\al'_f=\min_{x\ne 0}\{\al_{f,x}\}\q\h{with}\q\al_{f,x}=\min R_{f,x}\,(\les 1).$$
Let $\Ic_{0}$ be the ideal sheaf of $\{0\}\subset X=\C^n$, and $\Jc(X,\al D)$ be the multiplier ideal sheaf, see \cite{ELSV}, \cite{La}. By \cite{Mu}, \cite{Sa3}, we have
$$\Jc(X,\al D)=\Ic_{0}^{k}\,\,\,\,\,\h{with}\,\,\, k=[d\al]-n+1,\,\,\h{if}\,\,\,\al<\al'_f\,(\les 1).
\leqno(1.8.1)$$
This is due to \cite{Mu} in the case of hyperplane arrangements. It implies that $j/d$ is a jumping coefficient of $D$ at $0$ for $n\les j\les d\al'_f$. (Recall that the jumping coefficients are rational numbers $\al$ such that $\Jc(X,\al D)\ne\Jc(X,\al'D)$ for any $\al'<\al$.)
Note that $j/d$ is a jumping coefficient outside the origin if $j/d=\al'_f$.
\sk
We denote also by $V$ the induced $V$-filtration on $\Oc_X\subset\B_f$, see (1.1). Combining (1.8.1) with \cite{Bu}, \cite{BuSa1}, we get
$$V^{\al}\Oc_X=\Ic_{0}^{k}\,\,\,\,\,\h{with}\,\,\, k=\lceil d\al\rceil-n,\,\,\h{if}\,\,\,\al\les\al'_f\,(\les 1).
\leqno(1.8.2)$$
Here $\lceil\beta\rceil:=\min\{k\in\Z\mid k\ges\beta\}$.
(Note that $V^{\al}\Oc_X$ and $\Jc(X,\al D)$ are slightly differently indexed so that $V^{\al}\Oc_X=\Jc(X,\al D)$ if and only if $\al$ is not a jumping coefficient.)
\sk
Let $\JC_f$ be the set of jumping coefficients of the divisor $D$. We have by \cite{Ko}, \cite{ELSV} and \cite{Sa3} respectively
$$\aligned &\al_f=\min\JC_f=\min R_f,\\ &\JC_f\cap (0,1)\subset R_f\cap (0,1),\\ &\JC_f\cap (0,\al'_f)=R_f\cap (0,\al'_f).\endaligned
\leqno(1.8.3)$$
\sk
Let ${\rm Sp}(f)=\msum_{\al\in\Q}\,n_{f,\al}\,t^{\al}$ be the {\it Steenbrink spectrum} of $f$, see \cite{St2}, \cite{St3} (and \cite{BuSa2} for the case of hyperplane arrangements). We have by definition
$$n_{f,\al}:=\msum_i\,(-1)^i\dim\Gr_F^pH^{n-1-i}(F_{\!f,0},\C)_{\be(-\al)}\q\h{with}\q p:=[n-\al],$$
where $F$ is the Hodge filtration, and
$$\be(\beta):=\exp(2\pi i\,\beta)\q(\beta\in\Q).
\leqno(1.8.4)$$
We can also define ${\rm Sp}^i(f)=\msum_{\al\in\Q}\,n^i_{f,\al}\,t^{\al}$ by
$$n^i_{f,\al}:=\dim\Gr_F^pH^{n-1-i}(F_{\!f,0},\C)_{\be(-\al)}\q\h{with}\q p:=[n-\al].$$
Similarly we have the {\it pole order spectrum}
$${\rm Sp}_P(f)=\msum_{\al\in\Q}\,n_{P,f,\al}\,t^{\al}\q\h{together with}\q{\rm Sp}^i_P(f)=\msum_{\al\in\Q}\,n^i_{P,f,\al}\,t^{\al},$$
defined by replacing the Hodge filtration $F$ with the pole order filtration $P$.
\sk
By (1.8.2), the coefficients $n_{f,k/d}=n^1_{f,k/d}$ for $0<k/d<\al'_f$ are equal to $\binom{k-1}{n-1}$, that is,
$$n_{f,k/d}=\dim F^{n-1}H^{n-1}(F_{\!f,0},\C)_{\be(-k/d)}=\h{$\binom{k-1}{n-1}$}\,\,\,\,\,\h{if}\,\,\,\,\,0<k/d<\al'_f\,(\les 1).
\leqno(1.8.5)$$
Here $\binom{k-1}{n-1}=0$ if $k<n$. We have the inequality $n_{P,f,k/d}\les\binom{k-1}{n-1}$ by using
$$\aligned &\Omega_{\bP^{n-1}}^{n-1}\otimes_{\Oc}P_{0}\Lc^{(k/d)}\simeq\Oc_{\bP^{n-1}}(k-n),\\ &\dim\Gamma(\bP^{n-1},\Oc_{\bP^{n-1}}(k-n))=\h{$\binom{k-1}{n-1}$}.\endaligned
\leqno(1.8.6)$$
Combining this inequality with (1.8.5) and using the inclusion $F^{n-1}\subset P^{n-1}$, we get
$$n_{P,f,k/d}=\dim P^{n-1}H^{n-1}(F_{\!f,0},\C)_{\be(-k/d)}=\h{$\binom{k-1}{n-1}$}\,\,\,\,\,\h{if}\,\,\,\,\,0<k/d<\al'_f\,(\les 1).
\leqno(1.8.7)$$
\bs\bs
\vbox{\centerline{\bf 2. Cohomology of twisted de Rham complexes}
\bs\nin
In this section we review and partially generalize the theory of Aomoto complexes due to Esnault, Schechtman, Terao, Varchenko, Viehweg. This can be combined effectively with the construction in (1.5).}
\ms\nin
{\bf 2.1.~Twisted de Rham complexes.} Let $D$ be a central hyperplane arrangement in $X=\C^n$ ($n\ges 3$) with a reduced equation $f$ of degree $d$. Here central means that $D$ is the affine cone of a projective hyperplane arrangement $Z$ in $Y=\bP^{n-1}$. Note that the Bernstein-Sato polynomial of a global defining equation of the affine cone $D$ of a divisor on $\bP^{n-1}$ is equal to that of a local equation at $0\in D$, using the $\C^*$-action. Let $Z_i\,(1\les i\les d)$ be the irreducible components of $Z$ where $d=\deg Z$. By \cite{Br2}, \cite{ESV}, \cite{STV}, the cohomology of the local systems on $U:=Y\setminus Z$ in (1.5) can be calculated as follows:
\sk
Let $x_1,\dots,x_n$ be coordinates of $\C^n$ such that $Z_d=\{x_n=0\}$. Then the complement $Y'$ of $Z_d$ in $Y$ is identified with $\C^{n-1}$. Let $g_i$ be a polynomial of degree $1$ on $Y'=\C^{n-1}$ defining $Z'_i:=Z_i\cap Y'$. Put
$$\om_i=\ddd g_i/g_i\q\h{for}\q 1\les i\les d-1.$$
For $\al=(\al_1,\dots,\al_{d-1})\in\C^{d-1}$, set
$$g:=g_1\cdots g_{d-1},\q g^{\al}:=g_1^{\al_1}\cdots g_{d-1}^{\al_{d-1}},\q\om^{\al}:=\msum_{i<d}\,\al_i\om_i.$$
Let $\Oc_{Y'}g^{\al}$ be a free $\Oc_{Y'}$-module of rank $1$ on $Y'$ with formal generator $g^{\al}$. There is a regular singular integrable connection $\nabla$ such that for $\xi\in\Oc_{Y'}$
$$\nabla(\xi g^{\al})=(\ddd\xi)g^{\al}+\xi\om^{\al}g^{\al}.$$
Let $\Ac^{p}_{g,\al}$ be the $\C$-vector subspace of $\Gamma(U,\Omega_{U}^{p}\,g^{\al})$ generated by
$$\om_{i_1}\wdg\cdots\wdg\om_{i_{p}}\,g^{\al}\q\h{for any}\q i_1<\cdots<i_{p}.$$
Then $\Ac^{\ssb}_{g,\al}$ with differential given by $\om^{\al}{\wedge}\,$ is a subcomplex of $\Gamma(U,\Omega_{U}^{\ssb}\,g^{\al})$. Put
$$\al_d=-\msum_{i<d}\,\al_i.$$
By \cite{Br2}, \cite{ESV}, \cite{STV}, we have the canonical quasi-isomorphism
$$(\Ac^{\ssb}_{g,\al},\om^{\al}{\wedge})\simto\Gamma(U,\Omega_{U}^{\ssb}\,g^{\al}),
\leqno(2.1.1)$$
if the following condition holds for any {\it dense} edge $L$ of $Z$:
$$\al_{L}:=\msum_{Z_i\supset L}\,\al_i\notin\N\setminus\{0\}.
\leqno(2.1.2)$$
In the case of a constant local system (that is, $\al_i=0$), this is due to Brieskorn \cite{Br2}, and we have
$$\dim\Ac^i_{g,\al}=b_i(U).
\leqno(2.1.3)$$
Under a condition stronger than the above one, the quasi-isomorphism (2.1.1) is shown in \cite{ESV} as a solution of Aomoto's conjecture, and it is shown in \cite{STV} that it is enough to assume condition (2.1.2) only for {\it dense} edges. If $Z$ is {\it generic} (that is, if $Z$ is a divisor with normal crossings), then condition (2.1.2) is equivalent to $\al_i\notin\N\setminus\{0\}$ for any $i\in\{1,\dots,d\}$ (since the dense edges consist of the $Z_i$ in this case), and \cite{ESV} is sufficient in this case.
\ms\nin
{\bf 2.2.~Partial generalization.} In the above argument, we assume $\sum_{i=1}^{d}\al_i=0$, that is, the $\Oc_{Y}$-module with meromorphic connection is a trivial line bundle. It is easy to satisfy the condition (2.1.2), if we can use a nontrivial line bundle $E$ with a meromorphic connection. In this case we have
$$\msum_{i=1}^{d}\,\al_i=-k\q\h{if}\,\,\, E\simeq\Oc_{Y}(k)\,\,\,\h{with}\,\,\, k\ges 0,$$
and we may even assume $\al_i\les 0$ for any $i$.
\sk
Let $(\tY,\tZ)\to (Y,Z)$ be an embedded resolution obtained by blowing up along the proper transforms of certain edges of $Z$. Then the pull-back $(\tE,\nabla)$ of $(E,\nabla)$ is a logarithmic connection, and $H^j(U,\Omega_{U}^{\ssb}\,g^{\al})$ is calculated by the hyper\-cohomology of the logarithmic de Rham complex $\DR_{\log}(\tE)$ whose $p$th component is $\Omega_{\tY}^{p}(\log\tZ)\otimes_{\Oc}\tE$, see \cite{De1}. It is not clear if there is a simple formula as in (2.1) in this case. However, we can show the following for each $p$ by increasing induction on $n$:
$$H^j(\tY,\Omega_{\tY}^{p}(\log\tZ)\otimes_{\Oc}\tE)=0\q\h{if}\,\, j>0.
\leqno(2.2.1)$$
Note that it has nothing to do with the connection. If $k=0$ (that is, if $E=\Oc_{Y}$), then (2.2.1) follows from \cite{Br2} and the $E_1$-degeneration of the Hodge spectral sequence \cite{De2}, see \cite{ESV}. If $k>0$, take a sufficiently generic hyperplane $H$ of $Y$, and consider the pull-back to $\tY$ of the short exact sequence
$$0\to\Oc_{Y}(k-1)\to\Oc_{Y}(k)\to\Oc_{H}(k)\to 0.$$
This pull-back is exact, since $H$ intersects each edge of $Z$ transversally. This implies also that the pull-back $\tH$ of $H$ gives an embedded resolution of $(H,H\cap Z)$, and we get an exact sequence
$$0\to N^*_{\tH/\tY}\otimes_{\Oc}\Omega_{\tH}^{p-1}(\log\tZ_{H})\to\Omega_{\tY}^{p}(\log\tZ)\otimes_{\Oc}\Oc_{\tH}\to\Omega_{\tH}^{p}(\log\tZ_{H})\to 0,$$
where $\tZ_{H}=\tZ\cap\tH$ and $N^*_{\tH/\tY}$ is the conormal bundle of $\tH$ in $\tY$. The latter is isomorphic to the pull-back of $\Oc_{H}(-1)$ since $\tH$ is the total transform of $H$. So we can proceed by increasing induction on $n$. (The assertion is clear if $n=2$, that is, if $Y=\bP^1$, since $d>n$.)
\sk
By (2.2.1) the restriction to $U$ induces the canonical quasi-isomorphism
$$\Gamma(\tY,\Omega_{\tY}^{\ssb}(\log\tZ)\otimes_{\Oc}\tE)\simto\Gamma(U,\Omega_{U}^{\ssb}\,g^{\al}).
\leqno(2.2.2)$$
\bs\bs
\vbox{\centerline{\bf 3. $b$-Functions of hyperplane arrangements}
\bs\nin
In this section we first prove the main theorems, and then calculate the nearby and vanishing cycle sheaves in certain cases.}
\ms\nin
{\bf 3.1.~Proof of Theorem~1.} For the proof of $m_1=n$, note that $\om_{i_1}\wdg\cdots\wdg\om_{i_{n-1}}\ne 0$ for some $(i_1,\dots,i_{n-1})$ by hypothesis. Since this defines a nonzero logarithmic $(n-1)$-form of type $(n-1,n-1)$ on any embedded resolution of $(Y,Z)$, we get
$$\Gr^W_{2n-2}H^{n-1}(F_{\!f,0},\C)_1\ne 0,
\leqno(3.1.1)$$
and the assertion follows from (1.7) together with (1.1.3).
\sk
We now prove the assertion on $\max R_f$. We proceed by induction on $n$. We may assume that $D$ does not come from an arrangement in a lower dimensional vector space, and the assertion holds for the roots of $b_{f,x}(-s)$ at any $x\in D\setminus\{0\}$, taking a transversal space to each edge. Then by (1.4), it is enough to show
$$\Gr_{P}^jH^{n-1}(F_{\!f,0},\C)_{\be(-k/d)}=0\,\,\,\h{for}\,\,\begin{cases} j<n-2 &\h{if}\,\,\, k\in\{1,\dots,d-2\}\\ j<n-1 &\h{if}\,\,\, k\in\{d-1, d\},\end{cases}
\leqno(3.1.2)$$
where $\be(\beta)$ for $\beta\in\Q$ is as in (1.8.4). So the assertion for $k=d$ follows from the above argument. For each $k\in \{1,\dots,d-1\}$, we apply the argument in (2.2) to the case $\al_i=-k/d$ for any $1\les i\les d$ so that $E=\Oc_{Y}(k)$. Then (2.2.2) implies that
$$H^{n-1}(U,\Omega_{U}^{\ssb}g^{\-k/d})\,(=H^{n-1}(F_{\!f,0},\C)_{\be(-k/d)})$$
is generated by
$$\Gamma(\tY,\Omega_{\tY}^{n-1}(\log\tZ)\otimes_{\Oc}\tE),$$
which is identified with a subspace $\Gamma(Y,\Omega_{Y}^{n-1}\otimes_{\Oc}E(*Z))$. But this subspace is contained in $\Gamma(Y,\Omega_{Y}^{n-1}\otimes_{\Oc}P_1\Lc^{(k/d)})$ looking at the pole along the generic point of the proper transform of each irreducible component of $Z$. So the assertion follows if $k\in \{1,\dots,d-2\}$.
\sk
In the case $k=d-1$, set $\al_i:=1/d$ for $i\in \{2,\dots,d\}$, and $\al_1:=1/d-1$ so that $\sum_{i=1}^{d}\al_i=0$. Then (2.1.2) is satisfied (since $|m_{L}|<d$), and (2.1.1) holds. Let $V$ be the subspace of $\Ac_{g,\al}^{n-1}$ generated by
$$\om_{i_1}\wdg\cdots\wdg\om_{i_{n-1}}\,g^{\al}\q\h{for}\q\{i_1,\dots,i_{n-1}\}\subset\{2,\dots,d\}.$$
Then we have
$$\Ac_{g,\al}^{n-1}=\om_1\,{\wedge}\,\Ac_{g,\al}^{n-2}+V=\om^{\al}{\wedge}\,\Ac_{g,\al}^{n-2}+V,$$
and $V$ is a subspace of $\Gamma(Y,\Omega_{Y}^{n-1}\otimes_{\Oc}P_{0}\Lc^{(k/d)})$. So the assertion follows in this case. This finishes the proof of Theorem~1.
\ms\nin
{\bf 3.2.~Nearby cycles in the normal crossing case.} Let $f$ be a holomorphic function on a complex manifold $X$ such that $D:=f^{-1}(0)$ is a divisor with normal crossings. Let $m_i$ be the multiplicity of $f$ along each irreducible component $D_i$ of $D$. Let $\la\in\C^*$ of finite order $N>1$ (that is, $\la^{N}=1)$. Then for a general point $x$ of $D_{I}:=\bigcap_{i\in I}D_i$, it is well known that $H^j(F_{\!f,x},\C)_{\la}\ne 0$ if and only if $N$ divides $m_i$ for any $i\in I$, see \cite{St2}.
\sk
For $\la\in\C$, let $I(\la)=\{i\,|\,\la^{m_i}=1\}$. Let $W$ be the weight filtration on $\psi_{f,\la}\C_{X}$, which coincides with the monodromy filtration with center $n-1$, that is,
$$N^i:\Gr_{n-1+i}^W\psi_{f,\la}\C_{X}\simto\Gr_{n-1-i}^W\psi_{f,\la}\C_{X}\q\h{for}\,\, i>0,
\leqno(3.2.1)$$
where $N=\log T_{u}$ with $T_{u}$ the unipotent part of the monodromy $T$. (This holds without assuming the normal crossing condition.) Let $m_{x}$ be the smallest positive integer such that $N^{m_{x}}=0$ on the restriction of $\psi_{f,\la}\C_{X}$ to a sufficiently small neighborhood of $x$. Then
$$m_{x}=\#\{i\in I(\la)\,|\, x\in D_i\},
\leqno(3.2.2)$$
This follows from the construction of the weight filtration in \cite{St2} (see also \cite[Section 3.3]{Sa1}).
\ms\nin
{\bf 3.3.~Proof of Theorem~2.} We have an embedded resolution $(\tX,\tD)$ of $(X,D)$ by blowing up along the proper transforms of the dense edges of $D$ by increasing induction on the dimension of the edge as in \cite{STV} (choosing an order of the edges of the same dimension if necessary although the resolution does not depend on the order). Let $\tD_{L}$ be the proper transform of the exceptional divisor associated to the blowing-up along the proper transform of $L$. Here $\tD_{L}$ is the proper transform of $L$ if $L$ is an irreducible component of $D$. Since the multiplicity of the pull-back $\tf$ of $f$ along $\tD_{L}$ is $m_{L}$, the assertion follows from (3.2) and (1.1.3).
\ms
The following will be used in explicit calculations of Bernstein-Sato polynomials.
\ms\nin
{\bf 3.4.~Eigenvalues of the monodromy supported on the origin.} It is known that the Euler characteristic of $U:=\bP^{n-1}\setminus Z$ vanishes if and only if $D$ is {\it decomposable} (that is, there is a nontrivial decomposition $\C^n=\C^{n'}\times\C^{n''}$ such that $D$ is the union of the pull-backs of arrangements on $\C^{n'}$ and $\C^{n''}$), see \cite{STV}. This implies that, if the support of $\psi_{f,\la}\C_{X}$ is contained in $\{0\}$, then
$$H^{n-1}(F_{\!f,0},\C)_{\la}=0\,\,\Longleftrightarrow\,\,D\,\,\h{is decomposable.}
\leqno(3.4.1)$$
Here $\la=\be(-k/d)$ for some integer $k\in\{1,\dots,d-1\}$ (see (1.(.4)), and $\psi_{f,\la}\C_{X}$ is identified with $H^{n-1}(F_{\!f,0},\C)_{\la}$ where $H^j(F_{\!f,0},\C)_{\la}=0$ for $j\ne n-1$ since $\psi_{f,\la}\C_{X}$ is a shifted perverse sheaf supported on a point, see also \cite{CDO}. So (3.4.1) follows from the fact that the $H^j(F_{\!f,0},\C)_{\la}$ are calculated by the cohomology of a local system $L^{(k/d)}$ of rank $1$ on $U$ by (1.5) so that
$$\msum_j\,(-1)^j\dim H^j(F_{\!f,0},\C)_{\la}=\chi(U,L^{(k/d)})=\chi(U).
\leqno(3.4.2)$$
\ms\nin
{\bf 3.5.~Weight spectral sequence for vanishing cycles.} Assume $f$ is a homogeneous polynomial of $n$ variables with degree $d$, where $n\ges 3$. We denote by $\varphi_f\R_{h,X}[n{-}1]$ the mixed $\R$-Hodge module whose underlying $\R$-perverse sheaf is the vanishing cycle complex $\varphi_f\R_{X}[n{-}1]$, where $X:=\C^n$, see \cite{Sa1}. We give here a generalization of a construction in \cite{DiSa4} where the singularities of $Z$ are assumed to be ordinary double points.
\sk
For $\la\in\C^*$ with $|\la|=1$, we have by definition
$$\varphi_{f,\la}\C_{X}[n{-}1]={\rm Ker}(T_s-\la)\subset\varphi_{f}\C_{X}[n{-}1],$$
with $T_s$ the semisimple part of the monodromy $T$. (Here the kernel is taken in the abelian category of perverse sheaves.) We have the mixed $\R$-Hodge module
$$\varphi_{f,\la,\lab}\,\R_{X}[n{-}1],$$
whose underlying $\C$-perverse sheaf is given by
$$\varphi_{f,\la,\lab}\,\C_{X}[n{-}1]:=\begin{cases}\varphi_{f,\la}\C_{X}[n{-}1]\oplus\varphi_{f,\lab}\,\C_{X}[n{-}1]&\h{if}\,\,\,\la\ne\pm 1\\\varphi_{f,\la}\C_{X}[n{-}1]&\h{if}\,\,\,\la=\pm 1.\end{cases}$$
Similarly we can define $\tH^j(F_{\!f,0},\R)_{\la,\lab}$ which will be denoted by $\tH^j(F_{\!f,0})_{\la,\lab}$, where $\tH$ denotes the reduced cohomology, and $F_{\!f,0}$ is defined by $f^{-1}(1)$ in this case.
\sk
We have the weight spectral sequence of mixed $\R$-Hodge structures
$$_WE_1^{-i,i+j}=H^ji_0^*\Gr^W_i(\varphi_{f,\la,\lab}\,\R_{h,X}[n{-}1])\Longrightarrow\tH^{j+n-1}(F_{\!f,0})_{\la,\lab},
\leqno(3.5.1)$$
where $i_0:\{0\}\into X:=\C^n$ is the inclusion. This is induced by the weight filtration $W$ of the mixed $\R$-Hodge module $\varphi_{f,\la,\lab}\,\R_{h,X}[n{-}1]$. Here we use the canonical isomorphism of mixed $\R$-Hodge structures
$$H^ji_0^*(\varphi_f\R_{X}[n{-}1])=\tH^{j+n-1}(F_{\!f,0},\R),$$
compatible with the action of $T_s$. (This can be shown by using an argument similar to \cite[Section 1.6]{BuSa2}.)
\sk
Let $\Si_{\la}\subset Z$ such that
$$C(\Si_{\la})=\supp\varphi_{f,\la}\C_{X}[n{-}1]=\supp\varphi_{f,\lab}\,\C_{X}[n{-}1].$$
Here $C(V)\subset X$ denotes the affine cone of a closed subvariety $V\subset\bP^{n-1}$ in general. Assume
$$\dim\Si_{\la}=0.
\leqno({\rm C1})$$
We have the strict support decomposition of pure Hodge modules
$$\Gr^W_i(\varphi_{f,\la,\lab}\,\R_{h,X}[n{-}1])=\M^{\la,\lab}_{\{0\},i}\oplus\mopl_{z\in\Si_{\la}}\,\M^{\la,\lab}_{C(z),i}\q(i\in\Z).$$
Here $\M^{\la,\lab}_{V,i}$ is a pure $\R$-Hodge module of weight $i$ with strict support $V=\{0\}$ or $C(z)$ (where the latter means $C(\{z\})$). Since the weight filtration on the vanishing cycles is given by the shifted monodromy filtration (see \cite{Sa1}), we have the symmetry
$$\aligned N^i:\M^{\la,\lab}_{\{0\},r+i}&\simto \M^{\la,\lab}_{\{0\},r-i}(-i)\q(i\in\N),\\
N^i:\M^{\la,\lab}_{C(z),r+i}&\simto \M^{\la,\lab}_{C(z),r-i}(-i)\q(i\in\N),\endaligned
\leqno(3.5.2)$$
where the center of symmetry $r$ is given by
$$r:=\begin{cases}n-1&\h{if}\,\,\,\la\ne 1,\\ n&\h{if}\,\,\,\la=1.\end{cases}$$
\sk
For $z\in \Si_{\la}$, let $h_z$ be a local defining function of $(Z,z)$, and $F_{h_z}$ be the Milnor fiber of $h_z$. Assume the following condition:
$$\h{$h_z$ is a weighted homogeneous polynomial for any $z\in\Sigma_{\la}$.}
\leqno({\rm C2})$$
This implies that the monodromy on $H^{n-2}(F_{h_z},\C)_{\la}$ is just the multiplication by $\la$. We thus get $N=0$ on $\varphi_{f,\la}\C_X|_{X\setminus\{0\}}$, and the isomorphisms (3.5.2) implies
$$\M^{\la,\lab}_{C(z),i}=0\q\h{for}\,\,\,i\ne r,
\leqno(3.5.3)$$
since $\M^{\la,\lab}_{C(z),i}$ has strict support $C(z)$ (that is, its underlying perverse sheaf is an intersection complex with support $C(z)$). Moreover, the monodromy of the underlying local system of the restriction of $\M^{\la,\lab}_{C(z),r}$ to $C(z)\setminus\{0\}\cong\C^*$ is given by $T^{-d}$. (This is a well-known relation between the monodromy $T$ and the local system monodromy of $\psi_{f,\la}\C_X|_{C(z)\setminus\{0\}}$ for homogeneous polynomials $f$ of degree $d$. It can be verified, for instance, by using the point center blow-up at $0\in X$.) We then get
$$H^ji_0^*\M^{\la,\lab}_{C(z),r}=\begin{cases}H^{n-2}(F_{h_z})_{\la,\lab}&\h{if}\,\,\,j=-1\,\,\h{and}\,\,\,\la^d=1,\\0&\h{if}\,\,\,j\ne-1\,\,\,\h{or}\,\,\,\la^d\ne 1.\end{cases}
\leqno(3.5.4)$$
Note that $H^j(F_{f,0},\C)_{\la}=0$ unless $\la^d=1$. So we assume the following:
$$\la^d=1.
\leqno({\rm C3})$$
\sk
On the other hand we have
$$H^ji_0^*\M^{\la,\lab}_{\{0\},i}=\begin{cases}H^{\la,\lab}_{\{0\},i}&\h{if}\,\,\,j=0,\\0&\h{if}\,\,\,j\ne 0,\end{cases}$$
where the $H^{\la,\lab}_{\{0\},i}$ are pure $\R$-Hodge structures of weight $i$ such that
$$\M^{\la,\lab}_{\{0\},i}=(i_0)_*H^{\la,\lab}_{\{0\},i}.$$
These imply that $_WE_1^{-i,i+j}$ has weight $i+j$. In fact, it is shown in \cite{St2} that
$${\rm wt}\,H^{n-2}(F_{h_z})_{\la,\lab}=\begin{cases}n-2&\h{if}\,\,\,\la\ne1,\\n-1&\h{if}\,\,\,\la=1,\end{cases}$$
where ${\rm wt}\,H$ denotes the weight of a pure Hodge structure $H$.
\sk
The spectral sequence (3.5.1) then degenerates at $E_2$ (see \cite{De2}). Moreover
$$_W\ddd_1^{-i,i+j}:{}_WE_1^{-i,i+j}\to{}_WE_1^{1-i,i+j}\,\,\,\h{vanishes unless}\,\,\,(i,j)=(r,-1).
\leqno(3.5.5)$$
Here the only nonzero $E_1$-differential is
$$_W\ddd_1^{-r,r-1}:{}_WE_1^{-r,r-1}=\mopl_{z\in\Sigma_{\la}}\,H^{n-2}(F_{h_z})_{\la,\lab}\,\to {}_WE_1^{1-r,r-1}=H^{\la,\lab}_{\{0\},r-1}.$$
It is well-known that the monodromy $T$ on $\tH^j(F_{\!f,0},\C)$ is semisimple (since the geometric monodromy is given by $x_i\mapsto\zeta\,x_i$ with $\zeta:=\exp(2\pi i/d)$). Combining these, we then get
$$_W\ddd_1^{-r,r-1}\,\,\,\h{is surjective, and}\,\,\,H^{\la,\lab}_{\{0\},i}=0\,\,\,\h{for}\,\,\,|i-r|>1.
\leqno(3.5.6)$$
In fact, if $_WE_2^{1-r,r-1}\ne 0$, then $N\ne 0$ on $H^{n-2}(F_{\!f,0})_{\la,\lab}$ by (3.5.2) for $i=1$.
\sk
By (3.5.5) and (3.5.2) for $i=1$ we have the isomorphisms of mixed $\R$-Hodge structures
$$\Gr^W_{r+1}H^{n-1}(F_{\!f,0})_{\la,\lab}=H^{\la,\lab}_{\{0\},r+1}=H^{\la,\lab}_{\{0\},r-1}(-1).
\leqno(3.5.7)$$
So the surjectivity of $_W\ddd_1^{-r,r-1}$ together with the $E_2$-degeneration of the spectral sequence implies the following.
\ms\vbox{\nin
{\bf Proposition~3.6.} {\it In the notation of $(3.5)$ together with the hypotheses {\rm (C1--3)} in $(3.5)$, we have a short exact sequence of mixed $\R$-Hodge structures
$$0\to H^{n-2}(F_{\!f,0})_{\la,\lab}\to\mopl_{z\in\Sigma_{\la}}\,H^{n-2}(F_{h_z})_{\la,\lab}\to\Gr^W_{r+1}H^{n-1}(F_{\!f,0})_{\la,\lab}(1)\to 0,
\leqno(3.6.1)$$
compatible with the action of $T_s$.}}
\ms\nin
{\bf Remarks~3.7.} (i) Proposition~(3.6) is closely related with \cite[Theorem~1]{DiSa1}. In fact, the canonical self-pairing is non-degenerate on
$$\Gr^W_{n-1}H^{n-1}(F_{\!f,0})_{\la,\lab}=H_{\{0\},n-1}^{\la,\lab},$$
and the failure of the non-degeneration is given by the dual of $\Gr^W_nH^{n-1}(F_{\!f,0})_{\la,\lab}$.
\sk
(ii) The short exact sequence (3.6.1) implies, for instance, that we have $N\ne 0$ on $\varphi_{f,\la}\,\C_{X}$ if and only if
$$\dim H^{n-2}(F_{\!f,0},\C)_{\la}<\msum_{z\in\Sigma_{\la}}\,\dim H^{n-2}(F_{h_z},\C)_{\la}.
\leqno(3.7.1)$$
\ms
(iii) It may be conjectured that condition (3.7.1) is satisfied whenever $\Sigma_{\la}\ne\emptyset$. This is closely related to an argument in a recent preprint of R.~Kloosterman \cite{Kl}.
\ms
(iv) The relation between the constructions in (1.2) and (3.5) is not so simple. For instance, it does not seem trivial whether we have the equality
$$\Gr^{\al}_VM=(\D_X[\dd_tt])\,F_{-n}\Gr^{\al}_V\B_f\q\h{in}\,\,\,\Gr^{\al}_V\B_f.
\leqno(3.7.2)$$
If there is $u\in F_{-n}\Gr^{\beta}_V\B_f$ for $\beta<\al$ such that $(\dd_tt-\beta)^iu$ vanishes in $\Gr^{\beta}_V\B_f$, and is contained in $\Gr^{\al}_V\B_f$, then one would have to show that it is contained in the right-hand side of (3.7.2). This problem can be solved in the case of Theorem~(3.8) below by using Lemma~(3.9) below.
(There is a similar difficulty in the calculation of the induced Hodge filtration $F$ on $\Gr^{\al}_V\B_f$ in the normal crossing case \cite{Sa1}, and the argument is not so trivial as someone might imagine.)
\ms\vbox{\nin
{\bf Theorem~3.8.} {\it Under the assumptions~{\rm (C1--3)} in $(3.5)$, let $k\in\N$ with $\la:=\be(-k/d)\ne 1$ and $k/d\les\al'_f+1$, where $\al'_f$ is as in $(1.8)$. Then $m_{k/d}=2$ if we have for $q=[k/d]$}
$$\dim\Gr_F^{n-2-q}H^{n-2}(F_{\!f,0},\C)_{\la}<\msum_{z\in\Sigma_{\la}}\,\dim\Gr_F^{n-2-q}H^{n-2}(F_{h_z},\C)_{\la}.
\leqno(3.8.1)$$}
\ms\nin
{\it Proof.} We separate the proof into the three cases as follows:
\sk\nin
{\bf Case 1 : $k/d<1,\,\,q=0$.} It is enough to show the non-vanishing of the canonical morphism
$$\Gr_V^{k/d}M\to\Gr_n^W\Gr_V^{k/d}\B_f,
\leqno(3.8.3)$$
in the notation of (1.1). Here $W$ is the monodromy filtration with center $n-1$, which is associated to the operator $N=-(\dd_tt-k/d)$ as in (3.2.1). Note that we have $N^2=0$ on $\Gr_V^{k/d}\B_f$, since
$$\Gr_{n-1+i}^W\,\Gr_V^{k/d}\B_f=0\q\h{for}\,\,\,|i|>1,$$
by (3.5.3), (3.5.6). Here $r=n-1$ in (3.5) since $\la=\exp(-2pi k/d)\ne 1$.
\sk
The underlying filtered $\D$-module of the mixed Hodge module $\varphi_{f,\la,\lab}\,\R_{h,X}[n{-}1]$ is given by
$$\Gr_V^{k/d}(\B_f,F)\oplus\Gr_V^{1-k/d}(\B_f,F),$$
where the Hodge filtration $F$ on $\B_f$ is indexed like filtered right $\D$-modules so that
$$\min\{p\in\Z\mid F_p\B_f\ne 0\}=-n.$$
Since $F_{-n}\B_f=\Oc_X\subset M$, the non-vanishing of (3.8.3) is reduced to
$$F_{-n}\Gr^W_n\Gr_V^{k/d}\B_f\ne 0.
\leqno(3.8.4)$$
In the notation of (3.5), we have
$$\Gr^W_n\varphi_{f,\la,\lab}\,\R_{h,X}[n{-}1]=\M^{\la,\lab}_{\{0\},n}=(i_0)_*H^{\la,\lab}_{\{0\},n},$$
and (3.5.7) gives
$$\Gr^W_nH^{n-1}(F_{\!f,0})_{\la,\lab}=H^{\la,\lab}_{\{0\},n}.$$
Thus (3.8.4) is further reduced to
$$F^{n-1}\Gr^W_nH^{n-1}(F_{\!f,0},\C)_{\la}\ne 0,
\leqno(3.8.5)$$
(Note that $F$ is {\it shifted by} 1 when we take the nearby cycle functor $\psi_f$, that is, $\Gr_V^{\al}$ for $\al\in(0,1]$, see \cite{Sa1}.) 
So the assertion in the case $k/d<1$ follows from the short exact sequence in Proposition~(3.6).
\sk\nin
{\bf Case 2 : $k/d\in(1,\al'_f+1),\,\,q=1$.} Set $k':=k-d$. (Note that $k'<d$ since $\al'_f\les 1$.) In this case we have to determine $\Gr^{k'/d}_VM$, since (1.1.3) means that we have to know the order of nilpotency of the operator $N=-(\dd_tt-k/d)$ on
$$\Gr^{k/d}_VM\big/t(\Gr^{k'/d}_VM).$$
\sk
Since $f$ is a homogeneous polynomial, we have the Euler vector field $\xi$ such that $\xi f=f$, and hence
$$M=\D_Xf^s\q\h{in}\,\,\,\B_f.
\leqno(3.8.6)$$
By (1.8.2), we have
$$F_{-n}V^{k'/d}\B_f=V^{k'/d}\Oc_X=\Oc_X\,\C[x]_{k'-n}\q\h{for}\,\,\,\,k'/d\les\al'_f,$$
where $\C[x]_k$ denotes the purely degree $k$ part of the polynomial ring $\C[x]=\C[x_1,\dots,x_n]$. Using these, we can show the following by increasing induction on $k':$
$$\aligned V^{k'/d}M=\D_X\bl(\C[x]_{k'-n}\otimes 1\br)\subset V^{k'/d}\B_f\q\h{if}\,\,\,k'\les\al'_f\,d,\\ \Gr_V^{k'/d}M=\bl(\C[x]_{k'-n}\br)[\dd_1,\dots,\dd_n]\subset\Gr_V^{k'/d}\B_f\q\h{if}\,\,\,k'<\al'_f\,d,\endaligned
\leqno(3.8.7)$$
where $\dd_i:=\dd/\dd x_i$, and $\bl(\C[x]_{k'-n}\br)[\dd_1,\dots,\dd_n]$ can be viewed as the direct image of the vector space $\C[x]_{k'-n}$ under the inclusion $\{0\}\into X$ as a $\D$-module. In fact, (3.8.7) follows from Lemma~(3.9) below together with (3.8.6) by using the decreasing filtration $G$ on $\D_X$ and $M$ defined by
$$\aligned G^{k'}\D_X&=\begin{cases}\D_X&\h{if}\,\,\,k'\les n,\\ \D_X\bl(\C[x]_{k'-n}\br)&\h{if}\,\,\,n<k'\les\al'_f\,d,\\ 0&\h{if}\,\,\,k'>\al'_f\,d,\end{cases}\\ G^{k'}M&=V^{k'/d}M\q\h{for}\,\,\,k'\in\Z.\endaligned$$
Here the support of $\Gr_V^{k'/d}M$ is contained in the origin for $k'<\al'_f\,d$. This implies the second isomorphism in (3.8.7), which assures the injectivity assumption (3.9.1) in Lemma~(3.9) for the next step of the inductive argument on $k'$.
\sk
For $k'/d<\al'_f$, the argument in (3.5) implies moreover the inclusion
$$\Gr^{k'/d}_VM=\D_X(F_{-n}\Gr_V^{k'/d}\B_f)\subset W_{n-1}\,\Gr^{k'/d}_V\B_f,
\leqno(3.8.8)$$
together with the injectivity of the canonical morphism 
$$\Gr^{k'/d}_VM\into\Gr^W_{n-1}\,\Gr^{k'/d}_V\B_f,
\leqno(3.8.9)$$
where $W$ is the monodromy filtration on $\Gr^{k'/d}_V\B_f$ with center $n-1$. In fact, it is well-known that the Steenbrink exponents or spectral numbers (see \cite{St2}, \cite{St3}, and \cite{BuSa2}, \h{etc.}) of a weighted homogeneous polynomial coincide with the roots of $b_{h_z}(s)/(s+1)$ up to a sign by forgetting the multiplicities. (This follows from \cite{Ma1}, \cite[Chapters 1 and 4]{Sat}, \cite{ScSt}, \cite{St1}, etc.) So the condition $k'/d<\al'_f$ implies
$$\Gr_F^{n-2}H^{n-2}(F_{h_z},\C)_{\be(-k'/d)}=0,$$
and we get (3.8.8) by using (3.5.7) together with the surjectivity of the last morphism of (3.6.1). Then (3.8.9) follows, since $F_{-n}\Gr^W_{n-2}\,\Gr^{k'/d}_V\B_f=0$ by (3.5.2).
\sk
By (3.8.8--9), the quotient by $t(\Gr^{k'/d}_VM)$ does not affect the order of nilpotency of $N$ on $\Gr^{k/d}_VM$. So the assertion for $k'/d<\al'_f$ follows by an argument similar to the case $k<d$ and $q=0$.
\sk\nin
{\bf Case 3 : $k/d=\al'_f+1,\,\,q=1$.} Set $k':=k-d$ as above. In this case the last isomorphism of (3.8.7) and (3.8.8--9) do not hold. By the arguments in (3.5) we have
$$\Gr^W_{n-1+i}\Gr_V^{k'/d}(\B_f,F)=\begin{cases}\mopl_{z\in\Si_{\la}}\,\bl(M^{\la}_{C(z)},F\br)\oplus\bl(M^{\la}_{\{0\},n-1},F\br)&\h{if}\,\,\,i=0,\\ \bl(M^{\la}_{\{0\},n-1+i},F\br)&\h{if}\,\,\,i=\pm 1,\\ \,\,0&\h{if}\,\,\,|i|>1,
\end{cases}$$
where $\bl(M^{\la}_{\{0\},i},F\br)$, $\bl(M^{\la}_{C(z),i},F\br)$ are filtered regular holonomic $\D$-modules with strict support $\{0\}$ and $C(z)$ respectively, and $\la=\be(-k'/d)$. (For $C(z)$, see (3.5).)
\sk
By using an argument similar to the one in Cases~1 and 2, it is then enough to show
$$\Gr^W_{n-1+i}\,\Gr_V^{k'/d}M=\begin{cases}\mopl_{z\in\Si_{\la}}\,M^{\prime\,\la}_{C(z)}\oplus M^{\prime\,\la}_{\{0\},n-1}&\h{if}\,\,\,i=0,\\ M^{\prime\,\la}_{\{0\},n-1+i}&\h{if}\,\,\,i=\pm 1,\\ \,\,0&\h{if}\,\,\,|i|>1, 
\end{cases}
\leqno(3.8.10)$$
with
$$\aligned&M^{\prime\,\la}_{\{0\},n}=\bl(F_{-n}\Gr^W_n\Gr_V^{k'/d}\B_f\br)[\dd_1,\dots,\dd_n],\\&\q\q\q\,\,\,N:M^{\prime\,\la}_{\{0\},n}\simto M^{\prime\,\la}_{\{0\},n-2}.\endaligned
\leqno(3.8.11)$$
Here it is sufficient to show that there is a regular holonomic $\D_X$-submodule
$$M''\subset\Gr_V^{k'/d}\B_f,$$
such that $F_{-n}\Gr_V^{k'/d}\B_f\subset M''$, and (3.8.10--11) are satisfied with $\Gr_V^{k'/d}M$ replaced by $M''$.
(In fact, it is allowed to divide $\Gr^{k/d}_VM$ by a larger $\D_X$-submodule if we can show that the order of nilpotency is 2 by this.)
Note that the first condition is equivalent to the condition that $\Gr_V^{k'/d}M\subset M''$ by (3.8.7).
(Here we can also use the fact that the action of $s=-\dd_tt$ on $\Gr_V^{k'/d}\B_f$ is defined by using the Euler vector field $\xi$ with $\xi\,f=f$ as in (3.8.6) together with the commutation relation $[\xi,g]=(k''/d)g$ for a monomial $g$ of degree $k''=k'-n$. These imply that the $\D_X$-module constructed below coincides with $\Gr_V^{k'/d}M$.)
\sk
We have the filtered $\C$-vector spaces $\bl(H^{\la}_{\{0\},n-1+i},F\br)$ for $|i|\les 1$ such that
$$\bl(M^{\la}_{\{0\},n-1+i},F\br)=(i_0)_*\bl(H^{\la}_{\{0\},n-1+i},F\br)=\bl(H^{\la}_{\{0\},n-1+i}[\dd_1,\dots,\dd_n],F\br),$$
where $(i_0)_*$ is the direct image of filtered $\D$-modules, and the index of $F$ is defined like right $\D$-modules so that there is no shift of $F$ in the last term of the above equalities.
\sk
There are also filtered $\C$-vector spaces $\bl(H^{\la}_{C(z),n-1},F\br)$ for $z\in\Si_{\la}$ such that
$$\bl(M^{\la}_{C(z),n-1},F\br)=(i_{C(z)})_*\bl(\bl(H^{\la}_{C(z),n-1},F\br)\otimes_{\C}\Oc_{C(z)}\br).$$
Here $(i_{C(z)})_*$ is the direct image of filtered $\D$-modules under the inclusion $i_{C(z)}:C(z)\into X$, and $\bl(H^{\la}_{C(z),n-1},F\br)\otimes_{\C}\Oc_{C(z)}$ is the scaler extension of the filtered constant sheaf $(H^{\la}_{C(z),n-1},F\br)$ on $C(z)$, which can be viewed as a filtered regular holonomic $\D_{C(z)}$-module.
(The filtration $F$ on $H^{\la}_{C(z),n-1}\otimes_{\C}\Oc_{C(z)}$ is shifted by 1 from the standard one for right $\D$-modules.)
Note that the monodromy of $\bl(M^{\la}_{C(z),n-1},F\br)$ around $0\in C(z)$ is trivial (see a remark after (3.5.3)), and any polarizable variation of Hodge structure on $C(z)$ is constant.
\sk
By the short exact sequences associated with the weight filtration $W$, we have the extension classes
$$e_z\in{\rm Ext}_{\D_X}^1\bl(M^{\la}_{\{0\},n},M^{\la}_{C(z),n-1}\br),\q e'_z\in{\rm Ext}_{\D_X}^1\bl(M^{\la}_{C(z),n-1},M^{\la}_{\{0\},n-2}\br).$$
These can be identified with extension classes of $\D$-modules on $C(z)$.
They respectively correspond to
$$f_z\in{\rm Hom}_{\C}\bl(H^{\la}_{\{0\},n},H^{\la}_{C(z),n-1}\br),\q f'_z\in{\rm Hom}_{\C}\bl(H^{\la}_{C(z),n-1},H^{\la}_{\{0\},n-2}\br),$$
by using the morphisms
$$t:\Gr_V^0\to\Gr_V^1,\q\dd_t:\Gr_V^1\to\Gr_V^0,$$
for the extension $\D_{C(z)}$-modules corresponding to the extension classes $e_z$ and $e'_z$. Here $t$ is a coordinate of $C(z)$, and $V$ is the filtration of Kashiwara and Malgrange indexed by $\Z$ so that $\dd_tt-j$ is nilpotent on $\Gr_V^j$ for $j\in\Z$.
\sk
Note that, for any regular holonomic $\D_X$-modules $M_1,M_2$ supported at 0, we have
$${\rm Ext}_{\D_X}^j(M_1,M_2)=0\q\h{if}\,\,\,j\ne 0.
\leqno(3.8.12)$$
In particular, we have always
$$\msum_{z\in\Si_{\la}}\,e'_z\ssc e_z=0\q\h{in}\q{\rm Ext}_{\D_X}^j\bl(M^{\la}_{\{0\},n},M^{\la}_{\{0\},n-2}\br).$$
Note that this vanishing is essentially equivalent to the existence of an extension $\D_X$-module with graded quotients
$$M^{\la}_{\{0\},n-2},\q\mopl_z\,H^{\la}_{C(z),n-1}\oplus M^{\la}_{\{0\},n-1},\q M^{\la}_{\{0\},n},$$
where we use also (3.8.12) for $j=1$.
\sk
We see that $f_z$ preserves $F$ on $H^{\la}_{\{0\}}$, etc., and $f'_z$ preserves $F$ up to the shift by 1, since they are induced respectively by $t$ and $\dd_t$ in the above argument. Then the desired $\D_X$-submodule $M''\subset \Gr_V^{k'/d}\B_f$ can be constructed as an extension of
$$\aligned M^{\prime\prime\,\la}_{\{0\},n-1+i}&:=(i_0)_*\bl(F^nH^{\la}_{\{0\},n-1+i}\br)\subset M^{\la}_{\{0\},n-1+i}\q(i=0,1),\\ M^{\prime\prime\,\la}_{\{0\},n-2}&:=(i_0)_*\bl(F^{n-1}H^{\la}_{\{0\},n-2}\br)\subset M^{\la}_{\{0\},n-2},\\ M^{\prime\prime\,\la}_{C(z),n-1}&:=(i_{C(z)})_*\bl(F^nH^{\la}_{C(z),n-1}\otimes_{\C}\Oc_{C(z)}\br)\subset M^{\la}_{C(z),n-1},\endaligned$$
by restricting the above extension classes to these submodules (and using the theory of extensions of perverse sheaves in \cite{BBD}).
So the assertion follows in the case $k/d=\al'_f+1$.
This finishes the proof of Theorem~(3.8).
\ms
In the proof of the above proposition we used the following.
\ms\vbox{\nin
{\bf Lemma~3.9.} {\it Let $u:(A,G)\to(B,G)$ be a morphism of filtered abelian groups. Assume there are integers $a<b$ such that $G^aA=A$, $G^aB=B$, and we have the injectivity of
$$\Gr^p_Gu:\Gr^p_GA\to\Gr^p_GB\q\h{for}\,\,\,p<b.
\leqno(3.9.1)$$
Then}
$$u(A)\cap G^pB=u(G^pA)\q\h{for}\,\,\,p\les b.
\leqno(3.9.2)$$}
\ms\nin
{\it Proof.} Apply the snake lemma to the morphism between the short exact sequences
$$0\to G^{p+1}\to G^p\to \Gr_G^p\to 0,$$
for $A$, $B$. Then we get
$$u(G^pA)\cap G^{p+1}B=u(G^{p+1}A)\q\h{for}\,\,\,p<b.
\leqno(3.9.3)$$
So the assertion follows by increasing induction on $p$.
\ms\nin
{\bf Remarks~3.10.} (i) It is unclear whether the converse of Theorem~(3.8) holds, especially when the inequality (3.8.1) holds for some $q>[k/d]$. However, assuming conditions (C1--3) in (3.5), we have $m_{k/d}<2$ if the {\it equality} holds in (3.8.1) for any $q\ges[k/d]$. In fact, this is equivalent to the vanishing of certain graded pieces of $H^{\la}_{\{0\},n}$ in the proof of Theorem~(3.8).
\sk
(ii) Theorem~(3.8) is closely related to the following question of M.~Tomari: Does (0.3) hold only for generic central hyperplane arrangements? The answer to this question seems to be quite positive at least in the case $n=3$. This is closely related with Remark~(3.7)(iii) by Theorem~(3.8).
\bs\bs
\vbox{\centerline{\bf 4. Explicit calculations of Bernstein-Sato polynomials}
\bs\nin
In this section we show how to calculate the Bernstein-Sato polynomial in certain cases including the generic one.}
\ms\nin
{\bf 4.1.} Let $D$ be a central hyperplane arrangement in $X=\C^n$ with $n\ges 3$. This is the affine cone of a projective hyperplane arrangement $Z$ in $Y=\bP^{n-1}$. Assume $D$ is indecomposable so that $\chi(U)\ne 0$ with $U:=\bP^{n-1}\setminus Z$, see (3.4). Let $R_f$ be the set of the roots of $b_f(-s)$, and $R'_f$ be the union of the roots of $b_{f,x}(-s)$ for $x\in D\setminus\{0\}$. Put $\al'_f=\min R'_f$. Let $k\in\{1,\dots,d\}$. By (1.5.1) we have
$$\dim H^{n-1}(F_{\!f,0},\C)_{\be(-k/d)}=|\chi(U)|\q\h{if}\,\,\,\kod\notin R'_f+\Z,
\leqno(4.1.1)$$
where $\be(\beta)$ for $\beta\in\Q$ is as in (1.8.4).
Indeed, the last condition implies that the perverse sheaf $\psi_{f,\be(-k/d)}\C_X[n{-}1]$ is supported on $\{0\}$, and is identified with a vector space so that we have a canonical isomorphism
$$\psi_{f,\be(-k/d)}\C_X[n{-}1]=H^{n-1}(F_{\!f,0},\C)_{\be(-k/d)}\q\h{if}\,\,\,\kod\notin R'_f+\Z.
\leqno(4.1.2)$$
This implies that $H^j(F_{\!f,0},\C)_{\be(-k/d)}=0$ if $j\ne n-1$ and $\kod\notin R'_f+\Z$.
\sk
In the case $\kod\in [\al'_f,1)$, we {\it assume} in the notation of (2.1) that there is a subset $I$ of $\{1,\dots,d-1\}$ such that $|I|=k-1$ and condition (2.1.2) is satisfied for
$$\al_i=\begin{cases} 1-\kod\,\, &\h{if}\,\,\, i\in I\cup\{d\},\\-\kod &\h{if}\,\,\, i\in I^{c}:=\{1,\dots,d-1\}\setminus I,\end{cases}
\leqno(4.1.3)$$
so that (2.1.1) is a quasi-isomorphism. Note that the target of (2.1.1) calculates the $\la$-eigenspace of the Milnor cohomology $H^{\ssb}(F_{\!f,0},\C)_{\la}$ with $\la=\be(-k/d)$.
(There are many reasons for the definition of the $\al_i$ in (4.1.3), and it is not easy to modify this. For instance, the $\al_i$ should be greater than $-1$ for the relation with the pole order filtration in (1.5), and $\al_d$ should be positive when we consider the pole of $\om_i$ at infinity.)
\sk
Let $V(I)'$ be the vector subspace of $\Ac_{g,\al}^{n-1}$ generated by
$$\om_{i_1}\wdg\cdots\wdg\om_{i_{n-1}}\,g^{\al}\q\h{for}\,\,\,\{i_1,\dots,i_{n-1}\}\subset I,$$
and $V(I)$ be the image of $V(I)'$ in $H^{n-1}(\Ac_{g,\al}^{\ssb},\om^{\al}{\wedge})$, where $\al=(\al_i)$.
\sk
In the notation of (1.5), $V(I)'$ can be identified with a subspace of
$$\Gamma(Y,\Omega_{Y}^{n-1}\otimes_{\Oc}P_{0}\Lc^{(k/d)}).$$
This is closely related with a remark after (4.1.3).
\sk
We can determine whether $\kod$ and $\kod+1$ belong to $R_f$ in certain cases as follows.
\ms\vbox{\nin
{\bf Theorem~4.2.} {\it Let $k\in\{1,\dots,d\}$. With the above notation and assumption, we have the following\,$:$
\ms\nin
$(a)$ If $k=d-1$ or $d$, we have $1-\frac{1}{d},1\in R_f$ and $2-\frac{1}{d},2\notin R_f$.
\ms\nin
$(b)$ If $\kod<\al'_f$, then we have $\kod\in R_f$ if and only if $k\ges n$.
\ms\nin
$(c)$ If $\binom{k-1}{n-1}<\dim H^{n-1}(F_{\!f,0},\C)_{\be(-k/d)}$, then $\kod+1\in R_f$.
\ms\nin
$(d)$ If $\kod<\al'_f$, $\kod\notin R'_f+\Z$, and $\binom{k-1}{n-1}=|\chi(U)|$, then $\kod+1\notin R_f$.
\ms\nin
$(e)$ If $V(I)\ne 0$, then $\kod\in R_f$.
\ms\nin
$(f)$ If $V(I)=H^{n-1}(\Ac_{g,\al}^{\ssb},\om^{\al}{\wedge})$, then $\kod+1\notin R_f\setminus R'_f$.
\ms\nin
$(g)$ If $V(I)\ne H^{n-1}(\Ac_{g,\al}^{\ssb},\om^{\al}{\wedge})$ and $\dim V(I)'=\binom{k-1}{n-1}$, then $\kod+1\in R_f$.}}
\ms\nin
{\it Proof.} The assertions (a)--(d) follow from Theorem~1 together with (1.8) and (3.4). For the remaining assertions, note that by (1.6) and (2.1.1), $V(I)'$ and $V(I)$ are respectively identified with subspaces of
$$\Gamma(Y,\Omega_{Y}^{n-1}\otimes_{\Oc}P_{0}\Lc^{(k/d)}),\q P^{n-1}H^{n-1}(F_{\!f,0},\C)_{\be(-k/d)}.$$
Then the assertions follow from (1.4) and (1.8.6). This finishes the proof of Theorem~(4.2).
\ms\nin
{\bf 4.3.~The generic case.} In the case of generic central hyperplane arrangements, we have $\al'_f=1$, $R'_f=\{1\}$. So Walther's formula (0.3) follows from Theorems 1--2 together with Theorem~(4.2)(b), (c), (d), since we have in tis case (see \cite{CoSu}, \cite{OrRa})
$$\h{$|\chi(U)|=\binom{d-2}{n-1}$}.
\leqno(4.3.1)$$
\sk
The following will be used in the proof of later theorems.
\ms\vbox{\nin
{\bf Proposition~4.4.} {\it With the notation of {\rm (4.1)}, let $Z(I)=\bigcup_{i\in I}Z_i\subset Y=\bP^{n-1}$, where $I$ is as in {\rm (4.1.3)} so that $|I|=k-1$. Assume $k\ges n$ and $Z(I)\cup Z_d$ is a divisor with normal crossings on $\bP^{n-1}$. Then the last hypothesis of Theorem~{\rm (4.2)(g)} is satisfied, that is,}
$$\dim V(I)'=\h{$\binom{k-1}{n-1}$}.
\leqno(4.4.1)$$}
\ms\nin
{\it Proof.} The assumption implies
$$\mwdg_{i\in I\setminus I'}\,\om_i\ne 0\,\,\,\,\,\h{for any}\,\,\,\,I'\subset I\,\,\,\,\h{with}\,\,\,\,|I'|=|I|-n-1=k-n.$$
(This does not hold if we assume only that $Z(I)$ is a divisor with normal crossings on $\bP^{n-1}$.)
Note that $V(I)'$ is identified with a subspace of
$$\Gamma(Y,\Omega_{Y}^{n-1}\otimes_{\Oc}P_{0}\Lc^{(k/d)})\cong\Gamma(Y,{\Oc}_Y(k-n)),$$
where $Y=\bP^{n-1}$ (see (1.5.2) for the last isomorphism), and it is spanned by
$$\bl(\mwdg_{i\in I\setminus I'}\,\om_i\br)\,g^{\al}=c_{I'}\bl(\mprod_{i\in I'}\,g_i\br)\,\ddd y_1\wdg\dots\wdg\ddd y_{n-1}\,g^{\alt},$$
where $c_{I'}\in\C^*$, $y_1,\dots,y_{n-1}$ are the coordinates of $\C^{n-1}$, and $\alt=(\alt_i)\in\Q^{d-1}$ is defined by
$$\alt_i=-k/d\q(\forall\,i\in\{1,\dots,d-1\}).$$
So (4.4.1) is reduced to the following which is shown by increasing induction on $|I|\ges n$:
\ms\nin
(A)\,\,The $P_{I'}:=\prod_{i\in I'}g_i$ for $I'\subset I$ with $|I'|=|I|-n+1$ generate the vector space of polynomials of degree $|I|-n+1$ in $n-1$ variables.
\ms\nin
In fact, for each $i\in I$, the inductive hypothesis implies that the $P_{I'}$ for $I'\subset I\setminus\{i\}$ with $|I'|=|I|-n$ generate the vector space of polynomials of degree $\les |I|-n$ in $n-1$ variables, where the assertion is clear if $|I|=n$. This finishes the proof of Proposition~(4.4).
\ms\nin
{\bf 4.5.~Calculation in a non-generic case.} In the notation of (4.1), assume $n=3$ and
$$\mult_zZ\les 3\,\,\,(\forall\,z\in Z)\q\h{with}\q\mult_zZ=3\,\,\,(\exists\,z\in Z).
\leqno({\rm M3})$$
Under the assumption $n=3$, condition (M3) is equivalent in the notation of (4.1) to
$$R'_f=\bl\{\h{$\frac{2}{3}, 1,\frac{4}{3}$}\br\},\q\al'_f=\h{$\frac{2}{3}$}.$$
\sk
In the non-generic case, we have to use the calculations in (2.1) and (3.5) together with Theorems~(3.8) and (4.2), Propositions~(4.4) and (4.7--8) below in order to determine the roots $\al$ of $b_f(-s)$ and their multiplicities for $\al\ges\al'_f$. For the moment we can calculate the Bernstein-Sato polynomial by this method only in relatively simple cases.
\sk
Let $\nu_i$, $\nu'_i$ be the number of $i$-ple points of $Z$ and $Z':=Z\setminus Z_d$ respectively. We have $\nu_3\ne 0$, $\nu_i=\nu'_i=0$ ($i>3$) by condition~(M3).
It is well known (see \cite{OrSo}) that
$$b_0(U)=1,\q b_1(U)=d-1,\q b_2(U)=\nu'_2+2\nu'_3.
\leqno(4.5.1)$$
(The last equality is closely related to (4.5.3) below.)
We have moreover
$$\h{$\chi(U)=\binom{d-2}{2}-\nu_3$}.
\leqno(4.5.2)$$
In fact, if $\nu_3=0$, this is a special case of (4.3.1) (and can be shown by an elementary calculation of the Euler characteristic). The general case is reduced to this case by deforming slightly $Z$ and calculating the difference between the Euler numbers of
$$\{xy(x+y)=0\}\q\h{and}\q\{xy(x+y-1)=0\},$$
(see also the assertion in the end of \cite[Section~3.1]{BuSa2}).
\sk
In the notation of (2.1), if $g_i, g_j, g_{l}$ have a common zero, then
$$\om_i\wdg\om_j+\om_j\wdg\om_{l}+\om_{l}\wdg\om_i=0.
\leqno(4.5.3)$$
By \cite{OrSo} these are the only relations among the $\om_i\wdg\om_j$ for $i,j\in\{1,\dots,d-1\}$. (This is closely related to Lemma~(4.6) below.) Note that $\om_i\wdg\om_j=0$ if $Z'_i\cap Z'_j=\emptyset$ (that is, if $Z'_i, Z'_j$ are parallel in $\C^2$).
\sk
Let $\om^{\al}=\al_i\,\om_i+\al_j\,\om_j+\al_l\,\om_l$ $with$ $\al_i,\al_j,\al_l\in\C^*$. Then
$$\om^{\al}\wdg\om_i,\,\,\om^{\al}\wdg\om_j\,\,\,\h{are linearly independent}\iff\al_i+\al_j+\al_l\ne 0.
\leqno(4.5.4)$$
In fact, (4.5.3) implies
$$\om^{\al}\wdg\om_j-\om^{\al}\wdg\om_i=(\al_i+\al_j+\al_l)\,\om_i\wdg\om_j,
\leqno(4.5.5)$$
where $\om^{\al}\wdg\om_i$ and $\om_i\wdg\om_j$ are linearly independent, since so are $\om_i\wdg\om_j$ and $\om_j\wdg\om_l$ by the remark after (4.5.3).
\sk
Let $\om^{\beta}=\beta_i\,\om_i+\beta_j\,\om_j+\beta_l\,\om_l$ with $\beta_i,\beta_j,\beta_l\in\C$. If $\al_i+\al_j+\al_l\ne 0$, then (4.5.4) implies
$$\om^{\al}\wdg\om^{\beta}=0\,\,\Longrightarrow\,\,\om^{\beta}=c\,\om^{\al}\q\h{for some}\,\,c\in\C.
\leqno(4.5.6)$$
Indeed, setting $\eta'_i=\om^{\al}\wdg\om_i$, we get two relations
$$\al_i\,\eta'_i+\al_j\,\eta'_j+\al_l\,\eta'_l=0\q\h{and}\q\beta_i\,\eta'_i+\beta_j\,\eta'_j+\beta_l\,\eta'_l=0.$$
They must coincide up to a constant multiple since $\eta'_i,\eta'_j,\eta'_l$ span a subspace of dimension 2 by (4.5.4).
\ms
Note that the last assumption in Theorem~(4.2)(g) is satisfied only in the case where the arrangement is rather simple as in the case where the hypothesis of Proposition~(4.4) is satisfied. 
To verify the assumptions in Theorem~(4.2)(e), (f) and (g), we have to calculate
$$(\om^{\al}{\wedge}\,\Ac_{g,\al}^1)\cap V(I)'\q\h{or}\q(\om^{\al}{\wedge}\,\Ac_{g,\al}^1)+V(I)'\q\h{in}\q\Ac_{g,\al}^2.
\leqno(4.5.7)$$
For this we need some terminology as follows. Let $I$ be as in (4.1). Set
$$Z'(J):=\mcup_{i\in J}\,Z'_i\q\h{for}\,\,\,J\subset\{1,\dots,d-1\}.
\leqno(4.5.8)$$
\sk\vbox{\nin
{\bf Definitions.} (i) We say that $i,i'\in I$ are {\it strongly $I^c$-connected} if $Z'_i\cap Z'_{i'}\cap Z'(I^c)\ne\emptyset$.
\sk\nin
(ii) We say that $i,i'\in I$ are $I^c$-{\it connected} if there are $i_0,\dots,i_r$ ($r\ges 0$) such that $i_0=i$, $i_r=i'$, and $i_{l-1},i_l$ are strongly $I^c$-connected for any $l\in \{1,\dots,r\}$.
\sk\nin
(iii) We say that $I'\subset I$ is an {\it $I^c$-connected component} of $I$ if any $i,i'\in I'$ are $I^c$-connected and $I'$ is a maximal subset of $I$ satisfying this property.
\sk\nin
(vi) We say that an $I^c$-connected component $I'$ of $I$ is {\it contractible}, if the dual graph of $Z'(I')$ is contractible.
\sk\nin
(v) For an $I^c$-connected component $I'$ of $I$, we define $\sigma(I')$ to be the number of points of $Z'(I')\cap Z'(I^c)$ contained in the smooth part of $Z'(I')$.}
\sk
Here the dual graph of $Z'(I')$ for an $I^c$-connected component $I'$ of $I$ consists of the vertices corresponding to the elements of $I'$ together with the edges corresponding to the points of $\bl(\Sing Z'(I')\br)\cap Z'(I^c)$.
Note that $Z'_i\cap Z'_j$ may be non-empty even if $i$ and $j$ belong to different $I^c$-connected components of $I$. We will see that these notions are useful for the study of (4.5.7) in the relatively simple cases (assuming $n=3$ and condition~(M3) as above), see Propositions~(4.7) and (4.8) below.
(Note that non-contractible $I^c$-connected components do not appear in the case $d\les 7$ if we choose $Z_d$ and $I^c$ appropriately.)
\ms
The following is well-known in the theory of hyperplane arrangements (see \cite{OrSo}). We note here a short proof for the convenience of the reader who is not familiar with the theory.
\ms\vbox{\nin
{\bf Lemma~4.6.} {\it Assume $n=3$. Then, in the notation of {\rm (4.1)}, we have
$$b_0(U)=1,\q b_1(U)=d-1,\q b_2(U)=\msum_{z\in Z\setminus Z_d}\,(m_{Z,z}-1),$$
where $m_{Z,z}$ is the the number of the components $Z_i$ containing $z$.}}
\ms\nin
{\it Proof.} Since $\dim\Ac_{g,\al}^i=b_i(U)$ (see (2.1.3)), the assertion holds for $i\ne 2$. So it is enough to calculate $\chi(U)$. Since
$$\chi(U)=\chi(Y')-\chi(Z'),$$
(where $Z':=Y'\cap Z$) and $\chi(Y)=1$, it is enough to show
$$\chi(Z')=d-1-\msum_{z\in Z\setminus Z_d}\,(m_{Z,z}-1).$$
But this is easily verified by using the short exact sequence
$$0\to\Q_{Z'}\to\mopl_{i<d}\Q_{Z'_i}\to\mopl_{z\in Z'}\Q^{m_{Z,z}-1}\to 0,$$
since $\chi(Z'_i)=1$. This finishes the proof of Lemma~(4.6).
\ms\vbox{\nin
{\bf Proposition~4.7.} {\it Assume $n=3$ and {\rm (M3)} in $(4.5)$. Let $I$ be as in {\rm (4.1)} so that {\rm (2.1.2)} is satisfied for the $\al_i$ defined by {\rm (4.1.3)} where $k=d-2$. Assume $Z'(I^c)$ is connected, and $\al_i+\al_j+\al_{j'}\ne 0$ for any $i\in I$, $j,j'\in I^c$ with $Z'_i\cap Z'_j\cap Z'_{j'}\ne\emptyset$. Let $c$ be the number of $I^c$-connected components $I'$ of $I$ with $\sigma(I')=0$, see Definition~{\rm (v)} in $(4.5)$. Then
$$\dim\,(\om^{\al}{\wedge}\,\Ac_{g,\al}^1)\cap V(I)'\les c,\,\,\h{that is,}\,\,\dim V(I)\ges\dim V(I)'-c.
\leqno(4.7.1)$$
Here the equality holds if for any $I^c$-connected component $I'$ of $I$ with $\sigma(I')=0$, the following two conditions are satisfied: $I'$ is contractible, and $\al_i+\al_{i'}+\al_j\ne 0$ for some $i,i'\in I'$, $j\in I^c$ with $Z'_i\cap Z'_{i'}\cap Z'_j\ne\emptyset$.}}
\ms\nin
{\it Proof.} Let $\eta g^{\al}=\sum_{i<d}\beta_i\om_i \,g^{\al}\in\Ac_{g,\al}^1$, and assume
$$\om^{\al}{\wedge}\,\eta g^{\al}=\msum_{i,j}\,\gamma_{i,j}\om_i\wdg\om_j \,g^{\al}=\msum_{i,j\in I}\,\gamma'_{i,j}\om_i\wdg\om_j \,g^{\al}\in V(I)',
\leqno(4.7.2)$$
where $\gamma_{i,j}:=\al_i\beta_j-\al_j\beta_i\in\C$ for $i<j$ and $\gamma_{i,j}=\gamma'_{i,j}=0$ for $i\ges j$. We may have $\gamma'_{i,j}\ne\gamma_{i,j}$ for $i,j\in I$, and the assumption (4.7.2) is not equivalent to the condition that $\gamma_{i,j}=0$ for $\{i,j\}\not\subset I$, since there are some relations among the $\om_i\wdg\om_j$ as in (4.5.3).
\sk
If $Z'_i\cap Z'_j$ is a double point of $Z'$ and $\{i,j\}$ is not contained in $I$, then $\om_{i,j}\notin V(I)'$, and we get by (4.7.2) together with a remark after (4.5.3)
$$\gamma_{i,j}=\al_i\beta_j-\al_j\beta_i=0\q(i>j).
\leqno(4.7.3)$$
\sk
If $Z'_i\cap Z'_j\cap Z'_{j'}\ne\emptyset$ with $i\in I$, $j,j'\in I^c$, then there is $c\in\C$ (depending on $i,j,j'$) such that
$$\beta_i=c\,\al_i,\q\beta_j=c\,\al_j,\q\beta_{j'}=c\,\al_{j'},$$
by (4.5.6) and the hypothesis of the proposition. So we may assume
$$\beta_j=0\q\h{for any}\,\,\,j\in I^c,
\leqno(4.7.4)$$
replacing $\eta$ with $\eta-c\,\om^{\al}$ for some $c\in\C$, since $Z'(I^c)$ is connected. Here we get also
$$\beta_i=0\,\,\,\,\h{if}\,\,\,Z'_i\,\,(i\in I)\,\,\,\h{intersects}\,\,\,Z'(I^c)\,\,\,\h{at a smooth point of }\,\,Z'(I).
\leqno(4.7.5)$$
\sk
Let $I'$ be an $I^c$-connected component of $I$. Assume $i,i'\in I'$ are $I^c$-strongly connected, that is, $Z'_i\cap Z'_{i'}\cap Z'_j\ne\emptyset$ for some $j\in I^c$. Then we have by (4.5.3)
$$\om_i\wdg\om_j\,g^{\al}=\om_{i'}\wdg\om_j\,g^{\al}\,\,\,\,\h{\rm mod}\,\,\,\,\C\om_i\wdg\om_{i'}\,g^{\al}\,\,\bl(\subset V(I)'\br).$$
Hence
$$(\al_i\om_i+\al_{i'}\om_{i'}+\al_j\om_j)\wdg(\beta_i\om_i+\beta_{i'}\om_{i'})\,g^{\al}=\al_j(\beta_i+\beta_{i'})\,\om_j\wdg\om_i\,g^{\al}\,\,\,\,\h{\rm mod}\,\,\,\,\, V(I)',$$
where $\beta_j=0$ by (4.7.4). So we get by (4.7.2) together with a remark after (4.5.3)
$$\beta_i+\beta_{i'}=0.
\leqno(4.7.6)$$
Moreover, we have by (4.1.3)
$$\gamma'_{i,i'}=\al_i\beta_{i'}-\al_{i'}\beta_i-\al_j\beta_i=-(\al_i+\al_{i'}+\al_j)\beta_i.$$
This is non-zero in the case $\beta_i\ne 0$ and the last hypothesis in the proposition is satisfied. Note that $\gamma'_{i,i'}$ is well-defined if $Z_i\cap Z_{i'}\cap Z'(I^c)\ne\emptyset$ (by using the hypothesis~(M3)).
\sk
By (4.7.6) the $\beta_i$ for $i\in I'$ are at most uniquely determined by $\beta_{i_0}$ for any $i_0\in I'$, and they are really determined for any value of $\beta_{i_0}$ if $I'$ is a contractible $I^c$-connected component of $I$ with $\sigma(I')=0$.
On the other hand, if $I'$ is an $I^c$-connected component of $I$ with $\sigma(I')\ges 1$, then we have $\beta_i=0$ for any $i\in I'$ by (4.7.5--6). So the assertion follows. This finishes the proof of Proposition~(4.7).
\ms\vbox{\nin
{\bf Proposition~4.8.} {\it With the notation and the assumptions of {\rm (4.1)} and {\rm (4.5)}, assume $|I^c|=d-k=2$, any $I^c$-connected component $I'$ of $I$ is contractible, and has $\sigma(I')\les 1$, and moreover we have $\al_i+\al_{j_1}+\al_{j_2}\ne 0$ if $Z'_i\cap Z'_{j_1}\cap Z'_{j_2}\ne\emptyset$ for some $i\in I$, where $I^{c}=\{j_1,j_2\}$. Then $H^2(\Ac_{g,\al}^{\ssb},\om^{\al}{\wedge})=V(I)$ so that the hypothesis of Theorem~{\rm (4.2)(f)} is satisfied and hence we have $\kod+1\notin R_f\setminus R'_f$.}}
\ms\nin
{\it Proof.} We have to show
$$\Ac_{g,\al}^2=\om^{\al}{\wedge}\,\Ac_{g,\al}^1\,\,\,\h{mod}\,\,\, V(I)'.
\leqno(4.8.1)$$
We first show that we have for $i\in I$, $j\in I^c$
$$\om_j\wdg\om_i \,g^{\al}\in\om^{\al}{\wedge}\,\Ac^1_{g,\al}+V(I)'.
\leqno(4.8.2)$$
We will prove this for each $I^c$-connected component $I'$ of $I$ by using
$$\om^{\al}{\wedge}\,\om_i \,g^{\al}=\al_{j_1}\om_{j_1}\wdg\om_i \,g^{\al}+\al_{j_2}\om_{j_2}\wdg\om_i \,g^{\al}\,\,\,\h{mod}\,\,\, V(I)'\q(i\in I).
\leqno(4.8.3)$$
Since $|I^c|=2$ and condition~(M3) in (4.5) is assumed, $I'$ is linearly $I^c$-connected (since it is contractible by hypothesis). We say that $i$ is an {\it end} element of $I'$, if $Z'_i$ is parallel to $Z'_{j_a}$ for $a=1$ or $2$. We first consider the case where $I'$ has an end element (where $|I'|$ may be 1).
In this case we have $\om_{j_a}\wdg\om_i=0$, and (4.8.2) for $j=j_b$ with $b=3-a$ follows from (4.8.3).
Since $I'$ is linearly $I^c$-connected, we can show (4.8.2) inductively by using (4.8.3) for $i\in I'$.
\sk
Now we show (4.8.2) in the case where $I'$ has no end element. In this case we have $|I'|=1$, and $Z'_i\cap Z'_{j_1}\cap Z'_{j_2}\ne\emptyset$ for $\{i\}=I'$, since $I'$ is contractible and $\sigma(I')\les 1$ for any $I^c$-connected component $I'$ of $I$ by hypothesis. In this case (4.8.2) follows from (4.5.4) by calculating $\om^{\al}{\wedge}\,\om_{j_1}\,g^{\al}$ and $\om^{\al}{\wedge}\,\om_{j_2}\,g^{\al}$ as in (4.8.3), since (4.8.2) is shown for the other $I^c$-connected components $I'$ of $I$ having end elements by using (4.8.3) for $i\in I'$.
\sk
Finally it remains to show
$$\om_{j_1}\wdg\om_{j_2}\,g^{\al}\in\om^{\al}{\wedge}\,\Ac_{g,\al}^1+V(I)'.
\leqno(4.8.4)$$
In the case where $Z'_i\cap Z'_{j_1}\cap Z'_{j_2}\ne\emptyset$, the assertion follows from the above assertion by using (4.5.3). So it remains to treat the case where $Z'_{j_1}\cap Z'_{j_2}$ is a double point of $Z'$. In this case any $I^c$-connected component $I'$ of $I$ has an end element.
Then (4.8.4) follows by calculating $\om^{\al}{\wedge}\,\om_j \,g^{\al}$ or $\om^{\al}{\wedge}\,\om_{j'}\,g^{\al}$ as in the above case.
This finishes the proof of Proposition~(4.8).
\sk
Using the above propositions, we can now prove the following.
\ms\vbox{\nin
{\bf Theorem~4.9.} {\it Assume $n=3$, $d\ges 4$, and condition~{\rm (M3)} in $(4.5)$. Set $e:=\lceil 2d/3\rceil-1$, that is, $e<2d/3\les e+1$. Then we have a subset $J\subset\{3,\dots,2d-2\}$ such that
$$\{3,\dots,e\}\cup\{d-1,d,d+1,d+2\}\subset J,
\leqno(4.9.1)$$
\vskip-7mm
$$b_f(s)=(s+1)\,\prod_{i=2}^{4}\Bigl(s+\frac{i}{3}\Bigr)\,\prod_{j\in J}\Bigl(s+\frac{j}{d}\Bigr).
\leqno(4.9.2)$$
Here $J$ contains $\frac{2d}{3}$ and $\frac{4d}{3}$, if $\frac{d}{3}\in\N$ and $(3.7.1)$ for $\la=\be(\pm 2/3)$ is satisfied.}}
\ms\nin
{\it Proof.} If the last condition of the theorem is satisfied, then we have $m_{2/3}=m_{4/3}=2$ by Theorem~(3.8). In fact, we have in this case
$$\aligned F^0H^1(F_{\!f,x_j},\C)_{\be(2/3)}&=H^1(F_{\!f,x_j},\C)_{\be(-2/3)},\q F^1H^1(F_{\!f,x_j},\C)_{\be(2/3)}=0,\\ F^1H^1(F_{\!f,x_j},\C)_{\be(-2/3)}&=H^1(F_{\!f,x_j},\C)_{\be(-2/3)},\endaligned$$
so that (3.8.1) is reduced to (3.7.1), where $\al'_f=2/3$ by condition~(M3).
The assertion then follows from Theorem~1, Theorem~(4.2)(b) and (1.1.3). (Note that $d-1\in J$ since $2d-1\notin J$.) Here we use the fact that (3.7.1) is satisfied in the case $d=6$, see Remark~(4.10)(i) below. This finishes the proof of Theorem~(4.9).
\ms\nin
{\bf Remarks~4.10.} (i) In many cases $H^1(F_{\!f,0},\C)_{\be(-2/3)}$ is very small. We have
$$\dim H^1(F_{\!f,0},\C)_{\be(-2/3)}\les 1\q\h{for}\,\,\, d\les 8,
\leqno(4.10.1)$$
and Example~(4.13)(i) is the only example such that $\dim H^1(F_{\!f,0},\C)_{\be(-2/3)}\ne 0$ with $d\les 8$, see \cite{AB}, \cite{CoSu}, \cite{Di}, \cite{BDS} and also Example~(4.13)(iii) below. By (1.5.1) we have
$$H^j(F_{\!f,0},\C)_{\be(-2/3)}=0\,\,\,\,(j\in\N)\q\h{if}\q d/3\notin\N.
\leqno(4.10.2)$$
\ms
(ii) As far as calculated, we have always in Theorem~(4.9)
$$J=\{3,\dots,r_{\!f}\}\q\h{with}\,\,\,r_{\!f}=2d-3\,\,\,\h{or}\,\,\,2d-2.
\leqno(4.10.3)$$
It would be interesting whether this is true in general. Note that the corresponding assertion does not hold for jumping coefficients, see Example~(4.13)(ii) below.
\ms\vbox{\nin
{\bf Theorem~4.11.} {\it In the notation and the assumptions of Theorem~$(4.9)$, assume $d\les 8$. Then {\rm (4.9.2)} holds with $J=\{3,\dots,r_{\!f}\}$ where $r_{\!f}$ is either $2d-3$ or $2d-2$ as in $(4.10.3)$.}}
\ms\nin
{\it Proof.} We have to determine the subset $J\subset\{3,\dots,2d-2\}$ in Theorem~(4.9). Set $e:=\lceil 2d/3\rceil-1$ as above. By (4.9.1) it is enough to show
$$\kod\in R_f\,\,\,\h{if}\,\,\,k\in\{e+1,\dots,d-2\},\q\kod+1\in R_f\,\,\,\h{if}\,\,\,k\in\{3,\dots,d-3\}.
\leqno(4.11.1)$$
Calculating $e+1$ as in (4.11.3) below, we see that it is enough to consider the case
$$e+1=d-2\q\h{with}\q d=6,7,8.
\leqno(4.11.2)$$
\sk
By classifying affine line arrangements $Z'$ in $\C^2$ (adding lines one by one inductively), we can calculate the maximum of $\nu_3$ (together with $e$) for each $d\in\{4,\dots,8\}$ as follows:
$$\begin{array}{cllllllll}
d&:&4&5&6&7&8\\ e+1&:&3&4&4&5&6\\
\max\,\nu_3&:&1&2&4&6&7\end{array}
\leqno(4.11.3)$$
Here we assume that the line at infinity contains as many triple points as possible. In fact, for $d=4$, $Z'$ is essentially unique, and is defined by $h=x(x-1)y$. For $d=5$, there are two possibilities:
$$h=xy(x-1)(y-1)\,\,\,\,\h{with}\,\,\,\,\nu_3=2,\q\h{or}\q h=x(x-1)y\,\ell\,\,\,\,\h{with}\,\,\,\,\nu_3=1,$$
where $\ell$ is a general polynomial of degree 1.
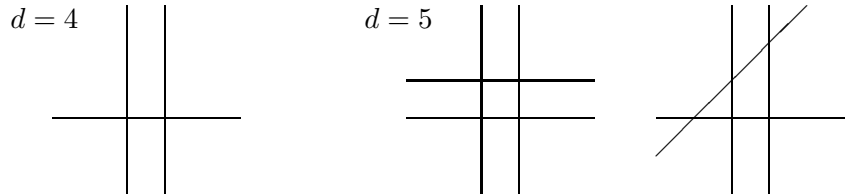
\sk
$$\raise2.2cm\h{\small$d=4$}\!\!\!\!\!
\h{$\setlength{\unitlength}{0.5cm}
\begin{picture}(5,5)
\put(0,2){\line(1,0){5}}
\put(2,0){\line(0,1){5}}
\put(3,0){\line(0,1){5}}
\end{picture}$}\q\q\q\q
\raise2.2cm\h{\small$d=5$}\!\!\!\!\!
\h{$\setlength{\unitlength}{0.5cm}
\begin{picture}(5,5)
\put(0,2){\line(1,0){5}}
\put(0,3){\line(1,0){5}}
\put(2,0){\line(0,1){5}}
\put(3,0){\line(0,1){5}}
\end{picture}$}\q\q
\h{$\setlength{\unitlength}{0.5cm}
\begin{picture}(5,5)
\put(0,2){\line(1,0){5}}
\put(2,0){\line(0,1){5}}
\put(3,0){\line(0,1){5}}
\put(0,1){\line(1,1){4}}
\end{picture}$}$$
\ms
We can continue like this for $d=6,7,8$. The details are left to the reader. These will be explained partially in later arguments.
\sk
From (4.5.2) we can deduce
$$\aligned\h{$\chi(U)-\binom{d-3}{2}$}&=d-3-\nu_3,\\ \h{$\chi(U)-\binom{d-4}{2}$}&=2d-7-\nu_3.\endaligned$$
Combined with (4.11.3) this implies
$$\h{$\chi(U)>\binom{k-1}{2}$}\q\h{if}\q\begin{cases}k=d-2\,\,\,\,\h{with}\,\,\,\,\nu_3<d-3,\,\,\,\h{or}\\ k<d-2\,\,\,\,\h{with}\,\,\,\,d>4.\end{cases}
\leqno(4.11.4)$$
By (1.5.1) we have furthermore
$$\dim H^2(F_{\!f,0},\C)_{\be(-k/d)}\ges\chi(U),
\leqno(4.11.5)$$
since $H^0(F_{\!f,0},\C)_{\be(-k/d)}=0$. Here the equality holds if and only if $H^1(F_{\!f,0},\C)_{\be(-k/d)}=0$. (The last condition is satisfied, if $d/3\notin\N$, see (4.1.2)).
\sk
By (4.11.4--5) for $k<d-2$ and Theorem~(4.2)(c), we get the last assertion of (4.11.1). For the first assertion of (4.11.1) we may assume
$$k=e+1=d-2\q\h{with}\q d=6,7,8.$$
In this case we have $c=0$ or $1$ in Proposition~(4.7), and it implies the first assertion of (4.11.1) by Theorem~(4.2)(e). This finishes the proof of Theorem~(4.11).
\ms\vbox{\nin
{\bf Theorem~4.12.} {\it In the notation and the assumption of Theorem~$(4.11)$, assume $d\les 7$. Then we have
$$r_{\!f}=2d-2\,\,\,\,\h{if}\,\,\,\,\nu_3<d-3,
\leqno(4.12.1)$$
and the converse holds except for the case in $(4.12.6)$ below where $d=7$, $\nu_3=4$ and $r_{\!f}=2d-2$.}}
\ms\nin
{\it Proof.} By Theorem~(4.11) it is enough to examine whether $2-2/d\in R_f$ or not. So we may assume
$$k=d-2.
\leqno(4.12.2)$$
If $\nu_3<d-3$, we have $2-\frac{2}{d}\in R_f$ by (4.11.4--5) for $k=d-2$ and Theorem~(4.2)(c). We may thus assume
$$\nu_3\ges d-3.
\leqno(4.12.3)$$
\sk
Assuming (4.12.2--3), we can verify the assertion as follows:
\sk\nin
{\bf Case 1 :} $d=4$ or $5$. In this case it is easy to apply Proposition~(4.8), and we get $r_{\!f}=2d-3$.
\sk\nin
{\bf Case 2 :} $d=6$. We can choose $I^c=\{j_1,j_2\}$ so that $Z_{j_1}\cup Z_{j_2}$ contains all the triple points of $Z$, and $Z_{j_1}\cap Z_{j_2}$ is a double point of $Z$. Then these conditions imply that (2.1.2) is satisfied for (4.1.3). Moreover, the assumption of Proposition~(4.8) is satisfied in case $\nu_3=3$ or $4$, and we get $r_{\!f}=2d-3$.
(Note that we have a strict inequality in (4.11.5) if $\nu_3=4$, see Example~(4.13)(i) below.)
\sk\nin
{\bf Case 3 :} $d=7$. The equivalent conditions of (4.5.4) is always satisfied in this case (since $d/3\notin\N)$.
\sk
If $\nu_3=6$, there is only one projective isomorphism class (see Example~(4.13)(ii) below), where we have $r_{\!f}=2d-3$ by Proposition~(4.8).
\sk
If $\nu_3=5$, there are two 1-parameter families (up to projective isomorphisms) defined by the following polynomials in $\C^2\subset\bP^2:$
$$\aligned &h=xy(x-1)(y-1)(x+y-1)(y-\la x)\q\h{($\la\in\C$ generic)},\\ &h=xy(x-1)(y-1)(\la x-y)((\la-1)x-(y-1))\q\h{($\la\in\C$ generic)}.\endaligned
\leqno(4.12.4)$$
\sk
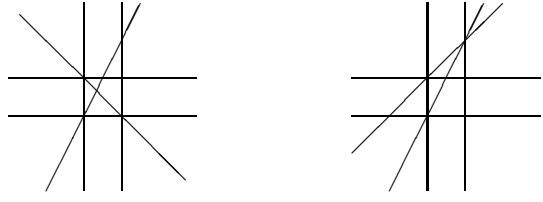
$$\h{$\setlength{\unitlength}{0.5cm}
\begin{picture}(5,5)
\put(0,2){\line(1,0){5}}
\put(0,3){\line(1,0){5}}
\put(2,0){\line(0,1){5}}
\put(3,0){\line(0,1){5}}
\put(0.3,4.7){\line(1,-1){4.4}}
\put(1,0){\line(1,2){2.5}}
\end{picture}$}\q\q\q\q\q
\h{$\setlength{\unitlength}{0.5cm}
\begin{picture}(5,5)
\put(0,2){\line(1,0){5}}
\put(0,3){\line(1,0){5}}
\put(2,0){\line(0,1){5}}
\put(3,0){\line(0,1){5}}
\put(1,0){\line(1,2){2.5}}
\put(0,1){\line(1,1){4}}
\end{picture}$}$$
\ms
This can be shown by counting the number of projective lines in $Z$ containing three triple points of $Z$. We can take $I^c$ corresponding to $\{xy=0\}$. Then we can apply Proposition~(4.8), and get $r_{\!f}=2d-3$.
\sk
If $\nu_3=4$, there are three 2-parameter families (up to projective isomorphisms) defined by the following polynomials in $\C^2\subset\bP^2:$
$$\aligned &h=xy(x-1)(y-1)(\la x-y)(\mu x-(y-1))\q\h{($\la,\mu\in\C$ generic)},\\ &h=xy(x-1)(y-1)(\la x-(y-1))((x-1)-\mu y)\q\h{($\la,\mu\in\C$ generic)},\endaligned
\leqno(4.12.5)$$
and
$$h=xy(x-1)(y-1)(x+y-1)(\la x-y+\mu)\q\h{($\la,\mu\in\C$ generic)}.
\leqno(4.12.6)$$
\sk
$$\h{$\setlength{\unitlength}{0.5cm}
\begin{picture}(5,5)
\put(0,2){\line(1,0){5}}
\put(0,3){\line(1,0){5}}
\put(2,0){\line(0,1){5}}
\put(3,0){\line(0,1){5}}
\put(0,4){\line(2,-1){5}}
\put(1,0){\line(1,2){2.5}}
\end{picture}$}\q\q\q\q\
\h{$\setlength{\unitlength}{0.5cm}
\begin{picture}(5,5)
\put(0,2){\line(1,0){5}}
\put(0,3){\line(1,0){5}}
\put(2,0){\line(0,1){5}}
\put(3,0){\line(0,1){5}}
\put(0,3.5){\line(2,-1){5}}
\put(0,1){\line(1,1){4}}
\end{picture}$}\q\q\q\q\
\h{$\setlength{\unitlength}{0.5cm}
\begin{picture}(5,5)
\put(0,2){\line(1,0){5}}
\put(0,3){\line(1,0){5}}
\put(2,0){\line(0,1){5}}
\put(3,0){\line(0,1){5}}
\put(0.3,4.7){\line(1,-1){4.4}}
\put(0,0.4){\line(3,2){5}}
\end{picture}$}$$
\ms
This can be shown by counting the number of projective lines in $Z$ containing three triple points of $Z$ and looking at the lines in $Z$ passing through two triple points of $Z$. We can take $I^c$ corresponding to $\{xy=0\}$. Then, for (4.12.5), we have $c=0$ and $12/7\notin R_f$, that is, $r_{\!f}=2d-3$, by Proposition~(4.8). For (4.12.6), however, we have $c=1$ in Proposition~(4.7), and hence $12/7\in R_f$, that is, $r_{\!f}=2d-2$, by Propositions~(4.4), (4.7), and Theorem~(4.2)(g), since the last condition of Proposition~(4.7) is satisfied in this case. This finishes the proof of Theorem~(4.12).
\ms\nin
{\bf Examples~4.13.} (i) In the notation of (2.1), let $n=3$, $d=6$, and
$$h=(x^2-1)(y^2-1)(x+y).
\leqno(4.13.1)$$
This is the simplest example with $H^{n-2}(F_{\!f,0},\C)_{\la}\ne 0$ for some $\la\ne 1$. In this case we have
$$\dim H^1(F_{\!f,0},\C)_{\be(\pm 1/3)}=1,\q\dim H^2(F_{\!f,0},\C)_{\be(\pm 1/3)}=3,$$
where $\nu_3=4$, $\chi(U)=2$.
We take $I^{c}$ corresponding to $\{(x+1)(y+1)=0\}\subset\C^2$ in the notation of Theorem~(4.2), where $k=d-2$. Then $Z(I)\cup Z_d\subset\bP^2$ has normal crossings so that $\dim V(I)'=3$ by Proposition~(4.4), and we have $c=1$ in Proposition~(4.7) for $k=4$. However, we have $\dim V(I)=3$ in this case. (Note that the last assumption in Proposition~(4.7) is not satisfied.) So we get $2-2/d\notin R_f$ by Theorem~(4.2)(f) (or Proposition~(4.8)).
\ms
(ii) In the notation of (2.1), let $n=3$, $d=7$, and
$$h=(x^2-1)(y^2-1)(x^2-y^2).
\leqno(4.13.2)$$
This is the only example with $d=7$ and $\nu_3=6$ up to a projective equivalence. In fact, it is projectively equivalent to the arrangement defined by the following polynomial in $\C^2\subset\bP^2:$
$$h=xy(x-1)(y-1)(x+y-1)(x+y-2).$$
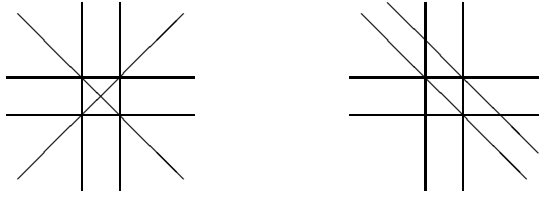
\sk
$$\h{$\setlength{\unitlength}{0.5cm}
\begin{picture}(5,5)
\put(0,2){\line(1,0){5}}
\put(0,3){\line(1,0){5}}
\put(2,0){\line(0,1){5}}
\put(3,0){\line(0,1){5}}
\put(0.3,4.7){\line(1,-1){4.4}}
\put(0.3,0.3){\line(1,1){4.4}}
\end{picture}$}\q\q\q\q\q
\h{$\setlength{\unitlength}{0.5cm}
\begin{picture}(5,5)
\put(0,2){\line(1,0){5}}
\put(0,3){\line(1,0){5}}
\put(2,0){\line(0,1){5}}
\put(3,0){\line(0,1){5}}
\put(0.3,4.7){\line(1,-1){4.4}}
\put(1,5){\line(1,-1){4}}
\end{picture}$}$$
\ms
In this case (4.9.2) holds with $J=\{3,\dots,11\}$ and $r_{\!f}=11=2d-3$. (Here we apply Proposition~(4.8) by setting $Z'(I^c)=\{(x+1)(y+1)=0\}$ for $k=d-2$ as in Example~(i).)
This coincides with a calculation by M.~Noro using his computer program Risa/Asir. In particular, $5/7$ is a root of $b_f(-s)$, although it is not a jumping coefficient, see \cite[Remark~(3.4)(ii)]{Sa3}, \cite[Theorem~3]{BuSa2}. Here we have $\nu'_3=4$ and $\chi(U)=4$.
\ms
(iii) In the notation of (2.1), let $n=3$, $d=9$, and
$$h=xy(y+2)(x-y)(x-y+1)(x+y-1)(x+y+2)(x-2y-1).
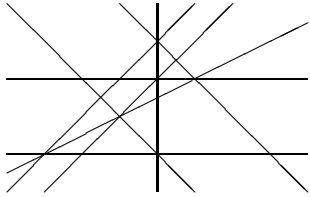
\leqno(4.13.3)$$
\sk
$$\h{$\setlength{\unitlength}{0.5cm}
\begin{picture}(8,5)
\put(0,3){\line(1,0){8}}
\put(4,0){\line(0,1){5}}
\put(0,1){\line(1,0){8}}
\put(1,0){\line(1,1){5}}
\put(0,0){\line(1,1){5}}
\put(3,5){\line(1,-1){5}}
\put(0,5){\line(1,-1){5}}
\put(0,0.5){\line(2,1){8}}
\end{picture}$}$$
\ms
This is the second simplest example with $\dim H^1(F_{\!f,0},\C)_{\be(\pm 1/3)}=1$, see \cite{AB}, \cite{CoSu}, \cite{Di}. In this case (4.9.2) holds with $J=\{3,\dots,r_{\!f}\}$. Here $r_{\!f}=15$ according to M.~Noro using the computer program Risa/Asir.)
\sk
For $k=6$, (4.4.1) holds, but (4.8.1) does not. Here $I^{c}$ corresponds to
$$\{y(x-y+1)(x+y+2)=0\}.$$
We have $b_1(U)=8$, $\dim V(I)'=10$, $\chi(U)=12$, and hence $0\ne V(I)\ne H^2(\Ac_{g,\al}^{\ssb},\om^{\al}{\wedge})$. So we get $6/9,15/9\in R_f$.
\sk
For $k=7$, it is not easy to show $16/9\notin R_f$ by using the method in this paper. However, this can be verified by using a different method, see Remark~(4.14)(iii) below.
\ms\nin
{\bf Remarks~4.14.} (i) The problem in Theorem~(4.2) is that the relation between
$$\Ac^{n-1}_{g,\al}\q\h{and}\q\Gamma(Y,\Omega_{Y}^{n-1}\otimes_{\Oc}P_{0}\Lc^{(k/d)})$$
may be rather complicated (both are subspaces of $\Gamma(Y,\Omega_{Y}^{n-1}\otimes_{\Oc}\Lc^{(k/d)}))$. If both (4.4.1) and (4.8.1) do not hold, we would have to enlarge the complex $(\Ac^{\ssb}_{g,\al},\om^{\al}{\wedge})$ as in (2.2). But it is very complicated to calculate this complex explicitly.
\sk
(ii) In the case where a hyperplane arrangement is sufficiently complicated, the hypotheses in Propositions~(4.4) and (4.7--8) are not satisfied, and we cannot calculate $b_f(s)$ by using the method in this paper. However, it is possible to apply a recent theory on pole order spectrum \cite{DiSa3} as is explained in Remark~(5.4)(iv) below so that the calculation of the roots of $b_f(s)$ can be reduced in certain cases to the one for the Hilbert series of the Jacobian ring.
Using this, we can calculate $b_f(s)$ for the $f$ in Examples~(4.13) and also in the proof of Theorem~4.12, see Remark~(iii) below.
We can also calculate $b_f(s)$ in the case of Walther's example \cite{Wa2} which showed that the Bernstein-Sato polynomials $b_f(s)$ are not combinatorial invariants of hyperplane arrangements, see Remark~(iv) below.
\sk
(iii) Let ${\rm Sp}^0_P(f)=\msum_kn^0_{P,f,k/n}\,t^{k/n}$ be the pole order spectrum for the highest Milnor cohomology, see (1.8). In the notation of \cite{DiSa3}, set
$$M:=H^n(\Omega^{\ssb},\ddd f\wdg),\q N:=H^{n-1}(\Omega^{\ssb},\ddd f\wdg),\q\mu_k=\dim M_k,\q\nu_k=\dim N_k,$$
with $(\Omega^{\ssb},\ddd f\wdg)$ the graded Koszul complex.
For the polynomials in (4.12.4), etc., we have the following (where $\la$, $\mu$ are specialized to some convenient integers):
\ms
\vbox{\nin For (4.12.4) with $\la=2$, we have $d=7$, $\nu_3=5$, $\nu_2=6$, $\tau=26$, $\chi(U)=5$, and
$$\begin{array}{rccccccccccccccccccccccccccccc}
k:&3& 4& 5& 6& 7& 8& 9& 10& 11& 12& 13& 14& \\
\mu_k:& 1 & 3 & 6 & 10 & 15 & 21 & 25 & 27 & 27 & 26 & 26 & 26 & \\
\nu_{k+7}:& & & 1& 5& 11& 16& 20& 23& 25& 26& 26& 26& \\
n^0_{P,f,k/n}:& 1 &3 &5 &5 &10 &5 &5 &4 & 2& & & & & & \\
\end{array}$$}
\sk
\vbox{\nin For (4.12.5) with $(\la,\mu)=(2,3)$, $(1,2)$, we have $d=7$, $\nu_3=4$, $\nu_2=9$, $\tau=25$, $\chi(U)=6$, and
$$\begin{array}{rccccccccccccccccccccccccccccc}
k:&3& 4& 5& 6& 7& 8& 9& 10& 11& 12& 13& 14& &\\
\mu_k:& 1 & 3 & 6 & 10 & 15 & 21 & 25 & 27 & 27 & 25 & 25 & 25 & \\
\nu_{k+7}:& & & & 4& 10& 15& 19& 22& 24& 25& 25& 25& \\
n^0_{P,f,k/n}:& 1 &3 &6 &6 &11 &6 &6 &5 & 3& & & \\
\end{array}$$}
\sk
\vbox{\nin For (4.12.6) with $(\la,\mu)=(-2,-1)$, we have $d=7$, $\nu_3=4$, $\nu_2=9$, $\tau=25$, $\chi(U)=6$, and
$$\begin{array}{rccccccccccccccccccccccccccccc}
k:&3& 4& 5& 6& 7& 8& 9& 10& 11& 12& 13& 14& \\
\mu_k:& 1 & 3 & 6 & 10 & 15 & 21 & 25 & 27 & 27 & 26 & 25 & 25 & \\ 
\nu_{k+7}:& & & & 4& 10& 15& 19& 22& 24& 25& 25& 25& \\
n^0_{P,f,k/n}:& 1 &3 &6 &6 &11 &6 &6 &5 & 3& 1& & \\
\end{array}$$}
\sk
\vbox{\nin For (4.13.1), we have $d=6$, $\nu_3=4$, $\nu_2=3$, $\tau=19$, $\chi(U)=2$, and
$$\begin{array}{rccccccccccccccccccccccccccccc}
k:&3& 4& 5& 6& 7& 8& 9& 10& 11& 12& \\
\mu_k:& 1 & 3 & 6 & 10 & 15 & 18 & 19 & 19 & 19 & 19 & \\
\nu_{k+6}:& & 1& 4& 9& 13& 16& 18& 19& 19& 19& \\
n^0_{P,f,k/n}:& 1 &3 &2 &6 &2 &3 &1 & & & & & \\
\end{array}$$}
\sk
\vbox{\nin For (4.13.2), we have $d=7$, $\nu_3=6$, $\nu_2=3$, $\tau=27$, $\chi(U)=4$, and
$$\begin{array}{rccccccccccccccccccccccccccccc}
k:&3& 4& 5& 6& 7& 8& 9& 10& 11& 12& 13& 14& \\
\mu_k:& 1 & 3 & 6 & 10 & 15 & 21 & 25 & 27 & 27 & 27 & 27 & 27 & \\
\nu_{k+7}:& & & 2& 6& 12& 17& 21& 24& 26 &27 &27 &27 \\
n^0_{P,f,k/n}:& 1 &3 &4 &4 &9 &4 &4 &3 & 1& & & & & & & \\
\end{array}$$}
\sk
\vbox{\nin For (4.13.3), we have $d=9$, $\nu_3=9$, $\nu_2=9$, $\tau=45$, $\chi(U)=12$, and
$$\begin{array}{rccccccccccccccccccccccccccccc}
k:&3& 4& 5& 6& 7& 8& 9& 10& 11& 12& 13& 14& 15& 16& 17& 18& \\
\mu_k:& 1 & 3 & 6 & 10 & 15 & 21 & 28 & 36 & 42 & 46 & 48 & 48 & 47 & 45 & 45 & 45& \\
\nu_{k+9}:& & & & 1& 3& 9& 17& 24& 30& 35& 39 &42 &44 &45 &45 &45 \\
n^0_{P,f,k/n}:& 1 &3 &6 &10 &12 &12 &19 &12 & 12& 12& 9& 6& 3& & & & & & \\
\end{array}$$}
\nin
Here the calculation of the $\mu_k$ is due to A.~Dimca and G.~Sticlaru.  (This can be done also by using Macaulay2, see \cite{Sa5}.)
We can easily calculate the $\nu_k$ from $\mu_k$, see \cite{DiSa3}.
These are compatible with the calculation of $b_f(s)$ in this paper by Remark~(5.4)(iv) below.
Note that $n^0_{P,f,(2d-2)/d}\ne 0$ and hence $(2d-2)/d\in R_f$ only in the case of (4.12.6).
\ms
(iv) Let $f$ be as in Walther's example in the degenerate case \cite{Wa2}, that is,
$$f=xyz(x+3z)(x+y+z)(x+2y+3z)(2x+y+z)(2x+3y+z)(2x+3y+4z).$$
We have $d=9$, $\nu_3=6$, $\nu_2=18$, $\tau=42$, $\chi(U)=15$, and
$$\begin{array}{rccccccccccccccccccccccccccccc}
k:&3& 4& 5& 6& 7& 8& 9& 10& 11& 12& 13& 14& 15& 16& 17& 18&\\
\mu_k:&1 &3 &6 &10 &15 &21 &28 &36 & 42& 46& 48& 48& 46& 43& 42& 42& \\
\nu_{k+9}:& & & & & 1& 6& 14& 21& 27& 32& 36& 39& 41& 42& 42& 42& \\
n^0_{P,f,k/n}:& 1 &3 &6 &10 &14 &15 &22 &15 & 15& 14& 12& 9& 5& 1& & \\
\end{array}$$
$$b_f(s)=(s+1)\,\mprod_{i=2}^4(s+i/3)\,\mprod_{i=3}^{16}(s+i/9).$$
\ms
Let $f$ be as in Walther's example in the non-degenerate case; for instance,
$$f=xyz(x+5z)(x+y+z)(x+3y+5z)(2x+y+z)(2x+3y+z)(2x+3y+4z).$$
We have $d=9$, $\nu_3=6$, $\nu_2=18$, $\tau=42$, $\chi(U)=15$ as above, and
$$\begin{array}{rccccccccccccccccccccccccccccc}
k:&3& 4& 5& 6& 7& 8& 9& 10& 11& 12& 13& 14& 15& 16& 17& 18& \\
\mu_k:&1 &3 &6 &10 &15 &21 &28 &36 & 42& 46& 48& 48& 46& 42& 42& 42& \\
\nu_{k+9}:& & & & & & 6& 14& 21& 27& 32& 36& 39& 41& 42& 42& 42& \\
n^0_{P,f,k/n}:& 1 &3 &6 &10 &15 &15 &22 &15 & 15& 14& 12& 9& 5& & & \\
\end{array}$$
$$b_f(s)=(s+1)\,\mprod_{i=2}^4(s+i/3)\,\mprod_{i=3}^{15}(s+i/9).$$
\ms\nin
The calculation of the $\mu_k$ of these two examples is also due to A.~Dimca and G.~Sticlaru.  (This can be done also by using Macaulay2, see \cite{Sa5}.)
We then get $b_f(s)$ by Remark~(5.4)(iv) below together with Theorem~(3.8). Here the $E_2$-degeneration of the pole order spectral sequence follows from \cite[Theorem~5.3]{DiSa3} and (1.8.2).
Note that $m_{15/9}=1$ by Remark~(3.10)(i).
\bs\bs
{\vbox{\centerline{\bf 5. Proof of Theorem~3.}
\bs\nin
In this section we terminate the proof of Theorem~3 by showing Proposition~1.}
\ms\nin
{\bf 5.1.~Proof of Proposition~1.} By Theorem~(3.8) it is enough to show (3.8.1) (with $q=0$) under the assumption
$$n=3,\,\,\,d=3m,\,\,\,k=3.
\leqno(5.1.1)$$
Here the right-hand side of (3.8.1) is equal to 1. (This follows, for instance, from the theory of spectra for hypersurface isolated singularities in the curve case, see \cite{St2}. It is also possible to use the multiplier ideal of $Z\subset\bP^2$ at a point of $Z$ with multiplicity $2m$, see (1.8).) The assertion is thus reduced to the following.
\ms\nin
{\bf Proposition~5.2.} {\it Under the assumptions of Proposition~$1$, we have}
$$H^1(F_{\!f,0},\C)_{\la}=0\q\h{for}\,\,\,\,\la=\exp(-2\pi i/m).
\leqno(5.2.1)$$
\ms\nin
{\it Proof.} This is proved by using the theory of Aomoto complexes as in (2.1). In this case we first show that there is a subset
$$I=\{i_1,i_2\}\subset\{1,\dots,d-1\},$$
such that condition (2.1.2) is satisfied by setting the $\al_i$ as in (4.1.3), that is
$$\al_i=\begin{cases} 1-\frac{1}{m}\,\, &\h{if}\,\,\, i\in I\cup\{d\},\\-\frac{1}{m} &\h{if}\,\,\, i\in I^{c}:=\{1,\dots,d-1\}\setminus I,\end{cases}$$
Here we assume that the point $z_0$ of $Z$ with multiplicity $2m$ is contained in the divisor at infinity (that is, $z_0\in Z_d=\bP^{n-1}\setminus\C^{n-1}$), and moreover $Z_{i_1}$ contains $z_0$. Set
$$J(z):=\bl\{j\in\{1,\dots,n-1\}\,\big|\,z\in Z_j\br\}.$$
If there is no point of $Z$ with multiplicity $m$, then $i_2$ may be any element in the complement of $J(z_0)$.
Assume there is a point $z_1\in Z$ with multiplicity $m$. If $J(z_0)\cap J(z_1)=\emptyset$, then $i_2$ may be any element of $J(z_1)$. In the other case, $\{i_2\}$ is the complement of $J(z_0)\cup J(z_1)$ in $\{1,\dots,d\}$. (In the last case we will assume $z_1\in Z_d$.)
\sk
We thus get the quasi-isomorphism (2.1.1), and it is enough to show
$$H^1(\Ac_{g,\al}^{\ssb},\om^{\al}{\wedge})=0.
\leqno(5.2.2)$$
By the theory on the kernel of the differential $\om^{\al}{\wedge}\,$ of the Aomoto complexes (see \cite{Fa}, \cite{LiYu}, \cite{FaYu}, etc.), the assertion is then reduced to the following (see \cite{BSY}, \cite[Section 1.5]{BDS}):
\ms\nin
{\bf Lemma~5.3.} {\it For any $i,j\in\{1,\dots,d-1\}$, there are $j_0,\dots,j_r$ such that $j_0=i$, $j_r=j$, and $Z'_{j_k}\cap Z'_{j_{k-1}}$ is a point of $Z'$ with multiplicity different from $m$ for any $k=1,\dots,r$, where $Z'_j=Z_j\setminus Z_d$.}
\ms\nin
{\it Proof.} The assertion is trivial if there is no point of $Z'$ with multiplicity $m$. Assume there is $z_1\in Z$ with multiplicity $m$.
Then we have either
$$J(z_0)\cap J(z_1)=\{d\}\q\h{or}\q J(z_0)\sqcup J(z_1)=\{1,\dots,d\},$$
and the assertion can be verified easily in both cases.
(For instance, $Z'=Z\setminus Z_d$ has no point of multiplicity $m$ in the first case.)
This finishes the proofs of Lemma~(5.3), Proposition~(5.2), Proposition~1, and Theorem~3.
\ms\nin
{\bf Remarks~5.4.} (i) It is show in \cite[Proposition~2.4]{BSY} that $-3/d$ is the only candidate for the pole of order $2$ of the topological zeta function $Z_{f,0}^{\rm top}(s)$ in the case of reduced central hyperplane arrangements with effective dimension $3$. Moreover it is really a pole of order 2 in this case if and only if there is a point of $Z$ with multiplicity $2d/3\in\N$.
\sk
(ii) It seems also possible to prove Proposition~(5.2) by using an argument in a recent preprint of R.~Kloosterman \cite{Kl} generalizing a result of \cite{Di} in the ordinary point case. This is also closely related with \cite{DiSa1}, \cite{DiSa3}.
\sk
(iii) Assume $f$ is a homogeneous polynomial and all the singularities of $Z$ are isolated and analytic-locally defined by weighted homogeneous polynomials. In this case we can show that $H''_{\!f,x}$ is $t$-torsion free, or equivalently, the pole order spectral sequence degenerates at $E_2$. (It will be proved in a forthcoming paper.) This gives another proof of \cite[Theorem~1.3]{BSY} for reduced hyperplane arrangements with $n=3$.
\sk
(iv) Under the assumption in Remark~(iii) above, the torsion-freeness of $H''_{\!f,x}$ implies that $\al$ is a root of $b_{f,0}(-s)$ if $n^0_{P,f,\al}\ne 0$, where $n^0_{P,f,\al}$ is as in (1.8). We can show conversely that if $\al$ is a root of $b_{f,0}(-s)$ such that $\al+i$ is not a root of $b_{f,x}(-s)$ for any $x\ne 0$ and $i\in\N$, then we have $n^0_{P,f,\al}\ne 0$. These follow from \cite[Theorem~2]{Sa3} (that is, Theorem~(1.4) in this paper).


\begin{thebibliography}{BuSa2}
\bibitem[AB]{AB} Artal-Bartolo, E., Combinatorics and topology of line arrangements in the complex projective plane, Proc. Amer. Math. Soc. 121 (1994), 385--390. 
\bibitem[BaSa]{BaSa} Barlet, D. and Saito, M., Brieskorn modules and Gauss-Manin systems for non-isolated hypersurface singularities, J.\ London Math.\ Soc.\ (2) 76 (2007), 211--224.
\bibitem[BBD]{BBD} Beilinson, A.A., Bernstein, J.N.\ and Deligne, P., Faisceaux pervers, Ast\'erisque 100, Soc. Math. France, Paris, 1982.
\bibitem[Be]{Be} Bernstein, J.N., The analytic continuation of generalized functions with respect to a parameter, Functional Analysis and its Applications 6 (1972), 273-285.
\bibitem[Br1]{Br1} Brieskorn, E., Die Monodromie der isolierten Singularit\"aten von Hyperfl\"achen, Manuscripta Math., 2 (1970), 103--161. 
\bibitem[Br2]{Br2} Brieskorn, E., Sur les groupes de tresses, S\'eminaire Bourbaki, 24\`eme ann\'ee (1971/1972), Exp. No. 401, Lect.\ Notes in Math. Vol.~317, Springer, Berlin, 1973, pp. 21--44.
\bibitem[Bu]{Bu} Budur, N., On Hodge spectrum and multiplier ideals, Math. Ann. 327 (2003), 257--270.
\bibitem[BDS]{BDS} Budur, N., Dimca, A.\ and Saito, M., First Milnor cohomology of hyperplane arrangements, in Topology of algebraic varieties and singularities, Contemp.\ Math., 538, Amer.\ Math.\ Soc., Providence, RI, 2011, 279--292.
\bibitem[BuSa1]{BuSa1} Budur, N. and Saito, M., Multiplier ideals, $V$-filtration, and spectrum, J. Alg.\ Geom.\ 14 (2005), 269--282.
\bibitem[BuSa2]{BuSa2} Budur, N. and Saito, M., Jumping coefficients and spectrum of a hyperplane arrangement, Math.\ Ann.\ 347 (2010), 545--579. 
\bibitem[BSY]{BSY} Budur, N., Saito, M.\ and Yuzvinsky, S., On the local zeta functions and the $b$-functions of certain hyperplane arrangements, with an appendix by W.\ Veys, J.\ London Math.\ Soc.\ (2) 84 (2011), 631--648. 
\bibitem[CDO]{CDO} Cohen, D.C., Dimca, A. and Orlik, P., Nonresonance conditions for arrangements, Ann. Inst. Fourier 53 (2003), 1883--1896. 
\bibitem[CoSu]{CoSu} Cohen, D. and Suciu, A., On Milnor fibrations of arrangements, J. London Math. Soc. 51 (1995), 105--119.
\bibitem[De1]{De1} Deligne, P., Equations Diff\'erentielles
\`a Points Singuliers R\'eguliers, Lect.\ Notes in Math. vol.~163, Springer, Berlin, 1970.
\bibitem[De2]{De2} Deligne, P., Th\'eorie de Hodge II, Publ. Math. IHES, 40 (1971), 5--58.
\bibitem[De3]{De3} Deligne, P., Le formalisme des cycles
\'evanescents, in SGA7 XIII and XIV, Lect.\ Notes in Math. 340, Springer, Berlin, 1973, pp. 82--115 and 116--164.
\bibitem[DeLo]{DeLo} Denef J.\ and Loeser, F., Caract\'eristiques d'Euler-Poincar\'e, fonctions z\'eta locales et modifications analytiques, J.\ Amer.\ Math.\ Soc.\ 5 (1992) 705--720. 
\bibitem[Di]{Di} Dimca, A., Singularities and Topology of Hypersurfaces, Springer, Berlin, 1992.
\bibitem[DiSa1]{DiSa1} Dimca, A.\ and Saito, M., Some consequences of perversity of vanishing cycles, Ann.\ Inst.\ Fourier 54 (2004), 1769--1792.
\bibitem[DiSa2]{DiSa2} Dimca, A.\ and Saito, M., A generalization of Griffiths' theorem on rational integrals, Duke Math.\ J.\ 135 (2006), 303--326.
\bibitem[DiSa3]{DiSa3} Dimca, A.\ and Saito, M., Koszul complexes and spectra of projective hypersurfaces with isolated singularities, arXiv:1212.1081.
\bibitem[DiSa4]{DiSa4} Dimca, A.\ and Saito, M., Generalization of theorems of Griffiths and Steenbrink to hypersurfaces with ordinary double points, arXiv:14034563.
\bibitem[ELSV]{ELSV} Ein, L., Lazarsfeld, R., Smith, K.E. and Varolin, D., Jumping coefficients of multiplier ideals, Duke Math.\ J.\ 123 (2004), 469--506.
\bibitem[ESV]{ESV} Esnault, H., Schechtman, V. and Viehweg, E., Cohomology of local systems on the complement of hyperplanes, Inv.\ Math.\ 109 (1992), 557--561.
\bibitem[Fa]{Fa} Falk, M., Arrangements and cohomology, Ann.\ Combin.\ 1 (1997), 135--157.
\bibitem[FaYu]{FaYu} Falk, M. and Yuzvinsky, S., Multinets, resonance varieties, and pencils of plane curves, Compos.\ Math.\ 143 (2007), 1069--1088.
\bibitem[Ka1]{Ka1} Kashiwara, M., $B$-functions and holonomic systems, Inv.\ Math.\ 38 (1976/77), 33--53.
\bibitem[Ka2]{Ka2} Kashiwara, M., Vanishing cycle sheaves and holonomic systems of differential equations, Lect.\ Notes in Math. 1016, Springer, Berlin, 1983, pp. 134--142.
\bibitem[Kl]{Kl} Kloosterman, R., On the relation between Alexander polynomials and Mordell-Weil ranks, equianalytic deformations and a variant of Nagata's conjecture (preprint).
\bibitem[Ko]{Ko} Koll\'ar, J., Singularities of pairs, Proc.\ Symp.\ Pure Math., A.M.S. 62 Part 1, (1997), 221--287.
\bibitem[La]{La} Lazarsfeld, R., Positivity in algebraic geometry II, Springer, Berlin, 2004.
\bibitem[LiYu]{LiYu} Libgober, A.\ and Yuzvinsky, S., Cohomology of the Orlik-Solomon algebras and local systems, Compos.\ Math.\ 121 (2000), 337--361.
\bibitem[Ma1]{Ma1} Malgrange, B., Le polyn\^ome de Bernstein d'une singularit\'e isol\'ee, Lect.\ Notes in Math.\ 459, Springer, Berlin, 1975, pp. 98--119.
\bibitem[Ma2]{Ma2} Malgrange, B., Polyn\^ome de Bernstein-Sato et cohomologie
\'evanescente, Analysis and topology on singular spaces, II, III (Luminy, 1981), Ast\'erisque 101--102 (1983), 243--267.
\bibitem[Mu]{Mu} Musta\c{t}\v{a}, M., Multiplier ideals of hyperplane arrangements, Trans.\ Amer.\ Math.\ Soc.\ 358 (2006), 5015--5023.
\bibitem[OrRa]{OrRa} Orlik, P.\ and Randell, R., The Milnor fiber of a generic arrangement, Ark. Mat. 31 (1993), 71--81.
\bibitem[OrSo]{OrSo} Orlik, P.\ and Solomon, L., Combinatorics and topology of complements of hyperplanes, Inv.\ Math.\ 56 (1980), 167--189.
\bibitem[Sa1]{Sa1} Saito, M., Mixed Hodge modules, Publ.\ RIMS, Kyoto Univ. 26 (1990), 221--333.
\bibitem[Sa2]{Sa2} Saito, M., On microlocal $b$-function, Bull.\ Soc.\ Math.\ France 122 (1994), 163--184.
\bibitem[Sa3]{Sa3} Saito, M., Multiplier ideals, $b$-function, and spectrum of a hypersurface singularity, Compos.\ Math.\ 143 (2007), 1050--1068.
\bibitem[Sa4]{Sa4} Saito, M., On real log canonical thresholds,
arXiv:0707.2308.
\bibitem[Sa5]{Sa5} Saito, M., {\it Hilbert series of graded Milnor algebras and roots of Bernstein-Sato polynomials}, arXiv:1509.06288.
\bibitem[Sat]{Sat} Sato, M. (ed.), Singularities of Hypersurfaces and $b$-Function (Proceedings of workshop in 1973), RIMS Kokyuroku 225 (in Japanese), 1975.
\bibitem[SatSh]{SatSh} Sato M.\ and Shintani, T., On zeta functions associated with prehomogeneous vector spaces, Proc.\ Nat.\ Acad.\ Sci.\ USA 69 (1972), 1081--1082.
\bibitem[STV]{STV} Schechtman, V., Terao, H. and Varchenko, A., Local systems over complements of hyperplanes and the Kac-Kazhdan conditions for singular vectors, J.\ Pure Appl.\ Algebra 100 (1995), 93--102.
\bibitem[ScSt]{ScSt} Scherk, J.\ and Steenbrink, J.H.M., On the mixed Hodge structure on the cohomology of the Milnor fibre, Math.\ Ann.\ 271 (1985), 641--665.
\bibitem[St1]{St1} Steenbrink, J.H.M., Intersection form for quasi-homogeneous singularities, Compositio Math.\ 34 (1977), 211--223.
\bibitem[St2]{St2} Steenbrink, J.H.M., Mixed Hodge structure on the vanishing cohomology, in Real and Complex Singularities (Proc. Nordic Summer School, Oslo, 1976) Alphen a/d Rijn: Sijthoff \& Noordhoff 1977, pp. 525--563.
\bibitem[St3]{St3} Steenbrink, J.H.M., The spectrum of hypersurface singularity, Ast\'erisque 179--180 (1989), 163--184.
\bibitem[Wa1]{Wa1} Walther, U., Bernstein-Sato polynomial versus cohomology of the Milnor fiber for generic hyperplane arrangements, Compos.\ Math.\ 141 (2005), 121--145.
\bibitem[Wa2]{Wa2} Walther, U., The Jacobian module, the Milnor fiber, and the $D$-module generated by $f^s$ (preprint).
\end{thebibliography}
\end{document}